\RequirePackage{fix-cm}
\documentclass[smallextended]{svjour3}
\smartqed
\usepackage{graphicx}
\usepackage{amsmath,amssymb,mathtools}
\usepackage[utf8]{inputenc}
\usepackage[T1]{fontenc}
\usepackage{enumitem}
\usepackage{etoolbox}
\AtEndEnvironment{proof}{\qed}

\newcommand{\bbc}{\mathbb{C}}
\newcommand{\bbn}{\mathbb{N}}

\newcommand{\bbr}{\mathbb{R}}
\newcommand{\bbe}{\mathbb{E}}
\newcommand{\bbp}{\mathbb{P}}
\newcommand{\bbF}{\mathbb{F}}
\newcommand{\cA}{\mathcal{A}}
\newcommand{\cL}{\mathcal{L}}
\newcommand{\cN}{\mathcal{N}}
\newcommand{\cJ}{\mathcal{J}}
\newcommand{\cI}{\mathcal{I}}
\newcommand{\cF}{\mathcal{F}}
\newcommand{\cR}{\mathcal{R}}
\newcommand{\cM}{\mathcal{M}}
\newcommand{\wt}{\widetilde}
\newcommand{\eps}{\varepsilon}

\journalname{~}
\renewcommand{\subclassname}{{\bfseries Mathematics Subject Classification (2020)}}

\begin{document}

\title{Reduction of the integro-differential equation to an
ordinary differential equation in ruin models
with investment in a risky asset}

\titlerunning{Reduction of the IDE to an ODE in ruin models with investment}

\author{Goluba Yurieva}

\authorrunning{G. Yurieva}

\institute{G. Yurieva \at
  Federal Research Center ``Computer Science and Control'' of the Russian Academy of Sciences, Moscow, Russia\\
  \email{goluba.yurieva@gmail.com}}

\date{Received: date / Accepted: date}

\maketitle

\begin{abstract}
We study ruin models with the surplus invested in a risky asset and Erlang-distributed claims: non-life insurance, annuity payments and the mixed model. Assuming that the integro-differential equation for the survival probability has a classical twice continuously differentiable solution, we prove that this solution satisfies a linear ordinary differential equation. We construct explicit annihilating operators and write the coefficients of the reduced equation in closed form through binomial coefficients (for the mixed model with equal shape parameters). A converse criterion is obtained: a bounded solution of the reduced equation satisfies the integro-differential equation if and only if the residual vanishes at zero initial capital together with its derivatives of order below that of the annihilator. The exponent of the power branch at infinity is determined by the parameters of the risky asset alone and does not depend on the shape parameter. For the mixed model we classify the fundamental set of solutions: the exponentially growing modes are excluded by boundedness. The connection coefficient of the power branch is nonzero by a known power-law lower bound for the ruin probability, and in the annuity model by the reduced equation and the boundary condition at zero. Hence the ruin probability has power-law asymptotics with a positive constant for all shape parameters and all positive exponents, and a complete asymptotic expansion up to a multiplicative constant in all three models. For non-life insurance with Erlang claims of order two the expansion coefficients are computed by an explicit recurrence.
\end{abstract}

\section{Introduction}

The analytical study of the ruin probability $\Psi(u)$ of an insurance
company that invests its surplus in a risky asset relies on the
integro-differential equation (IDE) for the survival function
$\Phi=1-\Psi$. The power of this approach depends on the
possibility of eliminating the nonlocal term and obtaining an ordinary
differential equation (ODE). For exponentially distributed jumps such
a reduction is known: in one-sided models it yields a third-order
equation for $\Phi$ (of second order for $\Phi'$)
\cite{FKP2002,KabPerg2016,PRYu2026}, and in the two-sided model, one of
fourth order \cite{KPuh2020}. In this paper we carry this reduction through for an
arbitrary Erlang distribution $\Gamma(k,\theta)$ and three types of
risk process. For an arbitrary $k$ the IDE implies the ODE. The
criterion in terms of the residual of the IDE
at zero covers the converse inclusion (Proposition~\ref{prop:converse}).
For bounded solutions of the ODE it reduces to at most $k$ conditions in
the non-life insurance model, holds automatically in the annuity model,
and requires at most $k_2$ conditions in the mixed model
(Remark~\ref{rem:converse}).

Interest in models with investments is due to the fact that they make
the tail of the ruin probability heavier: in contrast to the light tail
of the classical Cram\'er--Lundberg model, here $\Psi$ acquires a power
tail $\Psi(u)\sim Cu^{-\beta}$ with exponent $\beta=2a/\sigma^2-1$,
which is determined solely by the drift and the volatility of the asset
and does not depend on the claim distribution, provided that it has a
finite moment of order $\beta$ (for heavier claim tails the asymptotics
are determined by the claims themselves, see
\cite[Theorem~1]{P05143}). Discovered by Paulsen \cite{Paulsen1993} and
studied in detail in \cite{FKP2002,KalashnikovNorberg2002}, this effect
is investigated within two complementary approaches. The probabilistic
approach uses the Kesten--Goldie implicit renewal theory for the tails
of solutions of affine stochastic recursions
\cite{BDM2016,Goldie1991,Kesten1973} and its applications to insurance
with investments
\cite{EKS2022,KP2020,KP2023,KLP2026,Paulsen2002,PergamenshchikovZeitouni2006}:
it computes the exponent $\beta$ and, under the Cram\'er-type condition
on the logarithm of the price of the risky asset and moment conditions
on $P$ (both are satisfied for geometric Brownian motion and
Erlang-distributed jumps), the exact asymptotics
$\Psi(u)\sim Cu^{-\beta}$, $C>0$, for jumps of both signs
\cite{BKM2011,Nyrhinen2001,Paulsen2002}, while for a general L\'evy
return process only two-sided estimates of order $u^{-\beta}$ are
known \cite{KP2020}. The analytical approach, to which this paper belongs, treats
$\Phi$ as a solution of the IDE
\cite{FKP2002,Grandits2004,KPuh2020,PRYu2026} and aims at exact and
explicit results. For jump densities of both signs with a rational
Laplace transform (an Erlang distribution is among them), Antipov
\cite{Antipov2026} eliminated the integral term by a differential
operator in the variable $u$ and obtained an ODE of order $2n+2$ with
quadratic coefficients in $u$. Then, applying the Laplace transform and
a Tauberian theorem, he derived the asymptotics
$\Psi(u)\sim Cu^{-\beta}$, $C>0$, for $\beta\in(0,1)$.

The analytical approach requires smoothness of the survival function:
without it, writing the integro-differential equation in the classical
sense is unjustified. Existence of a classical $C^2$-smooth solution
under proportional investment was established by Promyslov for the annuity model
(life insurance) \cite{P2601} under
$\bbe[\xi]<\infty$, and for non-life insurance \cite{P05143} under the
minimal moment condition $\bbe[\xi^{\eps}]<\infty$ and continuity of
the distribution function. This strengthens the result of Grandits
\cite{Grandits2004} on the Sobolev regularity $W^{2,1}_{\mathrm{loc}}$
of the solution and the result \cite{AK2026} on $C^2$-smoothness under
$\bbe[\xi^{\beta}]<\infty$. Both papers \cite{P2601,P05143} also obtain the asymptotics
$\Psi(u)\sim Cu^{-\beta}$, $C>0$ (for non-life insurance, under
$\bbe[\xi^{\beta}]<\infty$). An alternative approach through viscosity
solutions is developed in \cite{P02518}. For exponentially distributed
jumps the problem is carried through to an exact solution in terms of
doubly confluent Heun functions in \cite{PRYu2026}.

Kabanov and Pukhlyakov \cite{KPuh2020} proved $C^2$-smoothness and
carried out such a reduction for a model with investments and two-sided
exponentially distributed jumps combining life and non-life insurance,
arriving at a linear fourth-order equation. The first reduction of the
IDE to an ODE for Erlang-distributed jumps appeared in the diploma thesis
of I.\,V.~Krechin, for annuity payments (see
\cite{Krechin2022}). In this paper, for the sake of completeness, we
provide a revised version of Krechin's reduction for this case. For
one-sided models with exponentially distributed jumps, an analogous
reduction is contained in \cite{FKP2002} (non-life insurance) and
\cite{KabPerg2016} (annuity payments), and its detailed exposition with
an exact solution through the doubly confluent Heun equation is given
in \cite{PRYu2026}. The method rests on a differential operator with
constant coefficients that annihilates the integral term: for
exponentially distributed one-sided jumps it is of first order (of second
order in the two-sided case), which leads to a third-order (respectively
fourth-order) equation whose asymptotics recovers the exponent
$\beta=2a/\sigma^2-1$. The possibility of such annihilation for all
distributions with a rational Laplace transform (Erlang, mixtures of
exponentials, phase-type distributions) dates back to the work of
Gerber and Shiu \cite{GerberShiu1998,GerberShiu2005} and was noted in
\cite{PRYu2026}. In full generality it was implemented by Antipov
\cite{Antipov2026}: for $k_1=k_2$ his equation (6) coincides with the
operator form \eqref{eq:ODE-C-op} of model~(C) (for arbitrary
$\theta_1,\theta_2$), and the reductions of one-sided models and the
case $k_1\neq k_2$ follow from the same construction. The contribution
of this paper with respect to \cite{Antipov2026} is as follows: closed
binomial formulas for the coefficients for one-sided models with an
arbitrary $k$ (for the mixed model with $k_1=k_2$ they are an explicit
specialization of equation (6) of \cite{Antipov2026} to the Erlang family), the operator form
for $k_1\neq k_2$, an asymptotic analysis by Thom\'e series for all
$\beta>0$ (without the restriction $\beta<1$) with a complete expansion
for non-life insurance, a converse criterion, and an analysis of the
singular point $u=0$.

The passage to the Erlang distribution $\Gamma(k,\theta)$ with an
arbitrary $k\ge 1$ does not reduce to replacing one density with
another and creates two new difficulties. First, the annihilating
operator is of order $k$, so the order of the resulting equation
grows together with $k$. Second, because of the investments the
coefficients of the equation depend on the capital and contain the
terms $u^2\Phi''$ and $u\Phi'$, so that an iterated application of an
operator of order $k$ via the Leibniz rule produces a large number of
cross terms, so a closed form of the coefficients requires exact
recurrence identities for the derivatives of the integral operator.

The main results of the paper are as follows.
\begin{enumerate}[label=(\roman*)]
  \item We construct annihilating differential operators $\cN$ and
  $\wt\cN$ of order $k$ with constant coefficients that reduce the
  integro-differential equation to a linear ordinary differential
  equation of order $k+2$ for one-sided models and $k_1+k_2+2$ for the
  mixed one.
  \item We write all coefficients of the reduced equations explicitly in closed
  form in terms of binomial coefficients and the parameters of the model
  (for the mixed model, in expanded form for
  $k_1=k_2$ and through the composition of operators in the general
  case), which makes the equations suitable for numerical and
  qualitative analysis for an arbitrary $k$.
  \item By the method of formal expansions near the irregular singular
  point $u=\infty$ (Thom\'e series) we prove that the exponent of the
  power branch does not depend on the shape parameter $k$ and coincides
  with the exponential case, $\beta=2a/\sigma^2-1$. The nonvanishing of
  the connection coefficient of this branch follows from the power-law
  lower bound for the ruin probability
  \cite[Proposition~4.1]{KPuh2020} (in the annuity model, from the ODE
  and the condition $\Phi(0+)=0$), so that
  $\Psi(u)\sim C_ku^{-\beta}$ with $C_k>0$ for all $k$ and all
  $\beta>0$ (Proposition~\ref{prop:realization}). For one-sided models
  this gives a different proof of the asymptotic part of the known
  results \cite{KP2020,P2601,P05143}. For the mixed model the asymptotics
  $\Psi(u)\sim Cu^{-\beta}$, $C>0$, are known from the probabilistic
  approach \cite{BKM2011,Nyrhinen2001,Paulsen2002}, while analytically
  it was obtained for $k_1=k_2=1$ \cite{KPuh2020} and for
  $\beta\in(0,1)$, $c\neq0$ \cite{Antipov2026}. The new elements are an
  analytical proof for all $k_1,k_2$, $\beta>0$, $c\in\bbr$ (apart from
  smoothness and the condition $\Phi(\infty)=1$, only the above lower
  bound is taken from probabilistic results) and a description of the
  solution space of the ODE. In all three models we obtain a complete
  asymptotic expansion up to a multiplicative constant $C_k$ (for
  non-life insurance it is known for $k=1$ from \cite{PRYu2026}, and
  for $k=2$ we give its coefficients by an explicit recurrence). For
  $k=1$ this agrees with \cite{FKP2002,PRYu2026} (in the annuity model,
  with \cite{KabPerg2016}).
  \item We obtain a converse criterion: a bounded solution of the
  reduced ODE satisfies the original IDE if and only if the residual of
  the IDE vanishes as $u\to0+$ together with derivatives up to order
  $K-1$, where $K$ is the order of the annihilator
  (Proposition~\ref{prop:converse}). The number of effective conditions
  is at most $k$ in model~(A), zero in model~(B), and at most $k_2$ in
  model~(C) (Remark~\ref{rem:converse}). For $k=1$ in model~(A) this is
  the known relation $c\Phi'(0+)=\alpha\Phi(0+)$ \cite{PRYu2026}.
\end{enumerate}
We conclude that the tail exponent is stable under deformation of
the jump-shape within the Erlang family: it is fixed by the parameters
of the risky asset $(a,\sigma)$.

The paper is organized as follows. In Sect.~\ref{sec:model} we
describe the model and write out integro-differential equations for the
survival function. In Sect.~\ref{sec:operators} we introduce
annihilating operators and establish the combinatorial identities
underlying them. In Sect.~\ref{sec:reduction} we prove reduction
theorems for all three models. In Sect.~\ref{sec:asymptotics} we
analyze formal solutions at the irregular singular point $u=\infty$ and
prove the realization of the power tail. Section~\ref{sec:examples}
contains the test cases $k=1$ and $k=2$, and
Sect.~\ref{sec:discussion} concludes with a discussion and open
questions.

\section{The model and assumptions}\label{sec:model}

All random objects are defined on a stochastic basis
$(\Omega,\cF,\bbF,\bbp)$ with a filtration $\bbF=(\cF_t)_{t\ge0}$
satisfying the usual conditions. It carries the following mutually independent objects:
a standard $\bbF$-adapted Wiener process $W=(W_t)$, a Poisson process
$N=(N_t)$ with intensity $\alpha>0$ and jump times
$0<T_1<T_2<\dots$, and an i.i.d.\ sequence $(\xi_i)_{i\ge1}$ of positive
random variables. In the mixed model we also have two independent
Poisson processes $N^{(1)}$, $N^{(2)}$ with intensities
$\alpha_1,\alpha_2>0$ (if $\alpha_1=0$ or $\alpha_2=0$ with the
corresponding sign of~$c$, the model degenerates into~(A) or~(B)) and
two sequences of jumps $\xi^{(1)}$, $\xi^{(2)}$.

The price of the risky asset $S$ is a geometric Brownian motion
\begin{equation*}
dS_t=S_t(a\,dt+\sigma\,dW_t),\qquad S_0>0,
\end{equation*}
with parameters $a\in\bbr$, $\sigma>0$. The return of the asset is
described by the process $dR_t=dS_t/S_t=a\,dt+\sigma\,dW_t$  The surplus
of a company investing the whole reserve in $S$ satisfies the
stochastic equation
\begin{equation*}
X_t=u+\int_0^t X_s\,dR_s+P_t,\qquad X_0=u\ge 0,
\end{equation*}
where $P=(P_t)$ is a risk process independent of $R$. The time of ruin
is $\tau^u:=\inf\{t\ge 0\colon X_t^u\le 0\}$, the survival function is
$\Phi(u):=\bbp(\tau^u=\infty)$, and the ruin probability is
$\Psi(u):=1-\Phi(u)$.

The risk process is one of the following three kinds.

\textbf{(A) Non-life insurance:}
\begin{equation}\label{eq:P-A}
P_t=ct-\sum_{i=1}^{N_t}\xi_i,\qquad c>0.
\end{equation}

\textbf{(B) Annuity payments:}
\begin{equation}\label{eq:P-B}
P_t=ct+\sum_{i=1}^{N_t}\xi_i,\qquad c<0.
\end{equation}

\textbf{(C) The mixed model:}
\begin{equation}\label{eq:P-C}
P_t=ct+\sum_{i=1}^{N_t^{(1)}}\xi_i^{(1)}-\sum_{i=1}^{N_t^{(2)}}\xi_i^{(2)},\qquad c\in\bbr,
\end{equation}
where $\xi_i^{(j)}\sim\Gamma(k_j,\theta_j)$, $j=1,2$, and the two
sequences are mutually independent.

We impose the following assumptions.

\begin{itemize}[leftmargin=2.4em]
  \item[\textbf{(A1)}] (\emph{non-degeneracy}) $a-\frac{1}{2}\sigma^2>0$,
  equivalently $\beta:=2a/\sigma^2-1>0$.
  \item[\textbf{(A2)}] (\emph{smoothness})
  $\Phi\in C^2((0,\infty))$.

  Existence of a classical solution for models~(A) and~(B) follows from
  independent results on non-life and annuity models:
  \cite[Theorem~1]{P05143} (under $\bbe[\xi^{\eps}]<\infty$ and
  continuity of the distribution function) and \cite[Theorem~1]{P2601}
  (under $\bbe[\xi]<\infty$). For the Erlang distribution both
  conditions hold.

  For the mixed model~(C), the smoothness of the survival function
  follows from \cite[Theorem~5.1]{KPuh2020}, and the validity of the IDE
  in the classical sense from \cite[Proposition~6.1]{KPuh2020}. The
  conditions on the jump densities and their derivatives are satisfied for
  the Erlang distribution. Indeed, the density
  \begin{equation}\label{eq:erlang}
  f_k(x)=\frac{x^{k-1}e^{-x/\theta}}{\theta^k(k-1)!},\qquad x\ge 0,
  \end{equation}
  and its derivatives $f_k'$, $f_k''$ are of the form ``a polynomial in
  $x$ multiplied by $e^{-x/\theta}$'', and therefore belong to
  $L^1(\mathbb{R}_+)$ for any $k\in\bbn$. The regularity is then
  bootstrapped within the paper itself, without circularity. By
  Lemma~\ref{lem:iter}, the highest derivative $J_k^{(k)}$ is expressed
  through the convolutions themselves $J_0=\Phi,\dots,J_k$ rather than
  through the derivatives of $\Phi$. Hence
  $\cJ_k\Phi,\cI_k\Phi\in C^{k+r}$ whenever $\Phi\in C^r$, and in
  particular $\cJ_k\Phi,\cI_k\Phi\in C^{k+2}$. On any compact subset of
  $(0,\infty)$, where $\frac{1}{2}\sigma^2u^2$ is bounded away from
  zero, the IDE of model~(A) gives
  $\Phi''=\bigl(\alpha\Phi-\alpha\cJ_k\Phi-(au+c)\Phi'\bigr)
  /(\frac{1}{2}\sigma^2u^2)\in C^1$, whence $\Phi\in C^3$. Repeating,
  we obtain $\Phi\in C^\infty((0,\infty))$ (in model~(B) similarly with
  $\cI_k\Phi$ and in model~(C) with $\cI_{k_1}\Phi$ and $\cJ_{k_2}\Phi$,
  and the boundedness of $\Phi$ guarantees the convergence
  of the integrals). The argument is local (smoothness is a local property), and
  therefore does not depend on the growth of the coefficients at
  infinity. It justifies the $k$-fold differentiation of the equations
  in Sect.~\ref{sec:reduction}. Furthermore, the coefficients of the
  reduced ODEs of Sect.~\ref{sec:reduction} are analytic and the leading
  coefficient does not vanish on $(0,\infty)$, and therefore $\Phi$ is
  real-analytic on $(0,\infty)$. Smoothness under more general
  conditions weakening $f',f''\in L^1$ (for distributions with a less
  regular density, for example) is the subject of a separate study and
  is not used in this paper.

  \item[\textbf{(A3)}] (\emph{non-degeneracy of the risk process}) the
  process $P$ is not an increasing process. In model~(A) this holds automatically
  ($\alpha>0$, downward jumps), in model~(B) because $c<0$, and in
  model~(C) because $\alpha_2>0$ (downward jumps). With unbounded
  (Erlang-distributed) downward jumps we have $\Psi(u)>0$ for all $u$.
  The case $\Psi\equiv0$ is possible only for an increasing $P$, cf.\
  \cite{KPuh2020}.
\end{itemize}

In what follows we use the following properties of the physical
survival function.
\begin{enumerate}[label=(\roman*)]
	\item $\Phi$ is nondecreasing: under a common noise
	$X^{u'}_t-X^{u}_t=(u'-u)\,\mathcal E_t>0$ for $u'>u$, where
	$\mathcal E_t=\exp\bigl((a-\frac{1}{2}\sigma^2)t+\sigma W_t\bigr)$.
	Hence $\tau^{u'}\ge\tau^{u}$, and in particular $G:=\Phi'\ge0$.
	\item $\Phi(\infty)=1$ under~(A1): in models~(A), (B) this is part
	of Theorem~1 in \cite{P05143} and \cite{P2601}, and in model~(C) of
	Lemma~2.1 in \cite{KPuh2020} (in the printed statement of that lemma
	the condition on the moment of order
	$\beta'\in(0,\beta\wedge1)$ of the downward claim size reads ``$<1$'',
	while only its finiteness is used in the proof). If the monetary unit
	is rescaled, $u,c,\xi\mapsto\lambda u,\lambda c,\lambda\xi$, the
	condition holds for a small $\lambda$, and the property
	$\Phi(\infty)=1$ is invariant: $\Phi_\lambda(\lambda u)=\Phi(u)$.
	\item $\Phi(0+)\in(0,1)$ in model~(A) \cite[Theorem~1]{P05143}, and
	$\Phi(0+)=0$ in model~(B) \cite[Theorem~1]{P2601}. In model~(C),
	$\Phi(0+)=0$ for $c<0$: before the first jump,
	$X_t=\mathcal E_t\bigl(u+c\int_0^t\mathcal E_s^{-1}ds\bigr)$ reaches
	zero in time of order $u/|c|$, while the probability of an upward
	jump before that moment tends to zero as $u\to0$ (the argument of
	\cite[Lemma~3]{P2601} carries over verbatim). For $c\ge0$ we have
	$\Phi(0+)>0$ (Proposition~\ref{prop:C-zero},
	Remark~\ref{rem:c0}).
	\item In model~(C) with $c\neq0$ the derivatives $\Phi'$, $\Phi''$
	are bounded on $(0,\infty)$ \cite[Theorem~5.1]{KPuh2020} (the
	statement of that theorem contains no proviso about $c$, but the
	boundedness in the proof is obtained precisely for $c\neq0$, and for
	$c=0$ it is in general false, see Remark~\ref{rem:c0}). For $c\neq0$
	it also follows from the local structure of solutions:
	Proposition~\ref{prop:C-zero} with the asymptotic existence theorem
	at zero, and Lemma~\ref{lem:actual} at infinity. For $c=0$ they are
	only locally bounded on $(0,\infty)$, and as $u\to0+$,
	$\Phi'=O(1+u^{\nu_+-1}|\ln u|)$ and
	$\Phi''=O(1+u^{\nu_+-2}|\ln u|)$ with the exponent $\nu_+>0$ from
	Remark~\ref{rem:c0}, which may be smaller than $1$ (the logarithmic
	factor is needed only for an integer $\nu_+$).
\end{enumerate}

Under assumptions~(A1)--(A3) the function $\Phi$ satisfies the
following IDEs (the derivation is standard, based on the generator of
the process, see \cite{KPuh2020,P05143}).

\textbf{Model~(A):} the boundary conditions $\Phi(z)=0$ for $z<0$ and
$\Phi(\infty)=1$, and
\begin{equation}\label{eq:IDE-A}
\frac{1}{2}\sigma^2 u^2\Phi''(u)+(au+c)\Phi'(u)+\alpha\!\int_0^u\!\Phi(u-y)f_k(y)\,dy-\alpha\Phi(u)=0,\quad u>0.
\end{equation}

\textbf{Model~(B):} the boundary conditions $\Phi(0)=0$ and
$\Phi(\infty)=1$, and
\begin{equation}\label{eq:IDE-B}
\frac{1}{2}\sigma^2 u^2\Phi''(u)+(au+c)\Phi'(u)+\alpha\!\int_0^\infty\!\bigl(\Phi(u+y)-\Phi(u)\bigr)f_k(y)\,dy=0.
\end{equation}

\textbf{Model~(C):} the boundary conditions $\Phi(z)=0$ for $z<0$ and
$\Phi(\infty)=1$, the behavior at zero being determined by the sign of
the drift (for $c>0$ survival from $u=0$ is possible, as in model~(A),
while for $c<0$ the drift takes the capital from zero into the negative
region and $\Phi(0+)=0$, as in~(B)), and
\begin{multline}\label{eq:IDE-C}
\frac{1}{2}\sigma^2 u^2\Phi''(u)+(au+c)\Phi'(u)
+\alpha_1\!\int_0^\infty\!\bigl(\Phi(u+y)-\Phi(u)\bigr)f_{k_1}^{(1)}(y)\,dy\\
+\alpha_2\!\int_0^u\!\Phi(u-y)f_{k_2}^{(2)}(y)\,dy-\alpha_2\Phi(u)=0.
\end{multline}

\section{Annihilating operators and basic identities}\label{sec:operators}

We introduce the second-order differential operator
\begin{equation}\label{eq:L}
\cL\Phi(u):=\frac{1}{2}\sigma^2 u^2\Phi''(u)+(au+c)\Phi'(u).
\end{equation}

For the Erlang density \eqref{eq:erlang} we define two integral
operators:
\begin{align}
\cJ_k\Phi(u)&:=\frac{1}{\theta^k(k-1)!}\int_0^u\Phi(u-y)y^{k-1}e^{-y/\theta}\,dy\notag\\
            &\phantom{:}=\frac{e^{-u/\theta}}{\theta^k(k-1)!}\int_0^u\Phi(z)(u-z)^{k-1}e^{z/\theta}\,dz,\label{eq:Jk}\\
\cI_k\Phi(u)&:=\frac{1}{\theta^k(k-1)!}\int_0^\infty\Phi(u+y)y^{k-1}e^{-y/\theta}\,dy\notag\\
            &\phantom{:}=\frac{e^{u/\theta}}{\theta^k(k-1)!}\int_u^\infty\Phi(z)(z-u)^{k-1}e^{-z/\theta}\,dz.\label{eq:Ik}
\end{align}
Notation: $J_k(u):=\cJ_k\Phi(u)$ and $I_k(u):=\cI_k\Phi(u)$, with the
convention $J_0:=\Phi$, $I_0:=\Phi$.

The annihilating differential operators of order $k$ with constant
coefficients are
\begin{equation}\label{eq:N}
\cN(f):=\sum_{n=0}^{k}\binom{k}{n}\theta^n f^{(n)},\qquad
\wt\cN(f):=\sum_{n=0}^{k}(-1)^n\binom{k}{n}\theta^n f^{(n)}.
\end{equation}
In symbolic form, $\cN=(I+\theta D)^k$ and $\wt\cN=(I-\theta D)^k$,
where $D=d/du$. In this notation the IDEs
\eqref{eq:IDE-A}--\eqref{eq:IDE-C} take the compact form
\begin{align}
\text{(A):}\quad &\cL\Phi+\alpha J_k-\alpha\Phi=0,\label{eq:compact-A}\\
\text{(B):}\quad &\cL\Phi+\alpha I_k-\alpha\Phi=0,\label{eq:compact-B}\\
\text{(C):}\quad &\cL\Phi+\alpha_1 I_{k_1}+\alpha_2 J_{k_2}-(\alpha_1+\alpha_2)\Phi=0.\label{eq:compact-C}
\end{align}

\begin{lemma}\label{lem:base}
For $u>0$ and $k\ge 1$:
\begin{equation}\label{eq:base-rules}
J_k'(u)=\frac{J_{k-1}(u)-J_k(u)}{\theta},\qquad
I_k'(u)=\frac{I_k(u)-I_{k-1}(u)}{\theta},
\end{equation}
with the convention $J_0=I_0=\Phi$.
\end{lemma}

\begin{proof}
We write $J_k$ in the equivalent representation
\[J_k(u)=\frac{e^{-u/\theta}}{\theta^k(k-1)!}\Lambda_k(u),\quad\text{where}\quad \Lambda_k(u)=\int_0^u\Phi(z)(u-z)^{k-1}e^{z/\theta}\,dz.\]
Differentiating, we get
\[
\Lambda_k'(u)=\Phi(u)(u-u)^{k-1}e^{u/\theta}+(k-1)\!\int_0^u\Phi(z)(u-z)^{k-2}e^{z/\theta}\,dz.
\]
For $k\ge 2$ the boundary term $\Phi(u)(u-u)^{k-1}e^{u/\theta}$ vanishes,
and the remaining integral equals
$(k-1)\theta^{k-1}(k-2)!\,e^{u/\theta}J_{k-1}(u)$. The case $k=1$ is
treated separately (bypassing the ambiguity $0^0$): then
$\Lambda_1(u)=\int_0^u\Phi(z)e^{z/\theta}\,dz$ and
$\Lambda_1'(u)=\Phi(u)e^{u/\theta}$ directly. In both cases
\[
J_k'(u)=-\frac{1}{\theta}J_k(u)+\frac{e^{-u/\theta}}{\theta^k(k-1)!}\Lambda_k'(u)=\frac{1}{\theta}(J_{k-1}-J_k),
\]
with the convention $J_0=\Phi$.

For $I_k$, by \eqref{eq:Ik}, we write
\[I_k(u)=\frac{e^{u/\theta}}{\theta^k(k-1)!}\cM_k(u),\]
where
\[\cM_k(u)=\int_u^\infty\Phi(z)(z-u)^{k-1}e^{-z/\theta}\,dz.\]
The integral converges because $\Phi$ is bounded, and the
differentiation under the integral sign is legitimate. For $k\ge2$ the
boundary term vanishes and
\begin{align*}
\cM_k'(u)&=-(k-1)\int_u^\infty\Phi(z)(z-u)^{k-2}e^{-z/\theta}\,dz\\
         &=-(k-1)\theta^{k-1}(k-2)!\,e^{-u/\theta}I_{k-1}(u),
\end{align*}
while for $k=1$, $\cM_1'(u)=-\Phi(u)e^{-u/\theta}$. In both cases
\[I_k'=\frac1\theta I_k+\frac{e^{u/\theta}}{\theta^k(k-1)!}\cM_k'=\frac1\theta(I_k-I_{k-1}).\]
\end{proof}

\begin{lemma}\label{lem:iter}
For $0\le n\le k$:
\begin{equation}\label{eq:iter-J}
J_k^{(n)}=\theta^{-n}\sum_{m=0}^{n}(-1)^{n+m}\binom{n}{m}J_{k-m},
\end{equation}
\begin{equation}\label{eq:iter-I}
I_k^{(n)}=\theta^{-n}\sum_{m=0}^{n}(-1)^{m}\binom{n}{m}I_{k-m}.
\end{equation}
\end{lemma}

\begin{proof}
Induction on $n$. The base $n=1$ is the content of
Lemma~\ref{lem:base}, rewritten as
$J_k'=\theta^{-1}\sum_{m=0}^{1}(-1)^{1+m}\binom{1}{m}J_{k-m}$ and
similarly for $I_k$.

\emph{Induction step for $J_k$.} Differentiating the assumed formula
for $n-1$ and applying Lemma~\ref{lem:base} to each $J_{k-m}$, we obtain
\[
J_k^{(n)}=\theta^{-n}\sum_{m=0}^{n-1}(-1)^{n-1+m}\binom{n-1}{m}(J_{k-m-1}-J_{k-m}).
\]
We split the sum into two:
\[
J_k^{(n)}=\theta^{-n}\Bigl[\sum_{m=0}^{n-1}(-1)^{n-1+m}\binom{n-1}{m}J_{k-m-1}-\sum_{m=0}^{n-1}(-1)^{n-1+m}\binom{n-1}{m}J_{k-m}\Bigr].
\]
In the first sum we change the index $m+1\mapsto m$. Both sums then
contain $J_{k-m}$ with the factors $(-1)^{n+m}\binom{n-1}{m-1}$ and
$(-1)^{n+m}\binom{n-1}{m}$, and the Pascal identity
$\binom{n-1}{m-1}+\binom{n-1}{m}=\binom{n}{m}$ for $1\le m\le n-1$ (the
extreme terms $m=0$ and $m=n$ enter once each, with
$\binom{n}{0}=\binom{n}{n}=1$) yields \eqref{eq:iter-J}. For $I_k$ the
proof is analogous with the replacement of signs:
$(-1)^{n+m}\mapsto(-1)^m$.
\end{proof}

\begin{remark}\label{rem:formal}
Below, we apply the operator $\cN=(I+\theta D)^k$ (and its mirror
$(I-\theta D)^k$) to functions of class $C^k$ on the half-line. To
formal series at the irregular singular point $u=\infty$ (Thom\'e
series, Sect.~\ref{sec:asymptotics}) we apply it termwise, and nowhere
below do we use the convergence of these series as an established
fact.
\end{remark}

\begin{lemma}\label{lem:annihilation}
For all $k\ge 1$:
\begin{equation}\label{eq:annih}
\cN J_k=\Phi,\qquad \wt\cN I_k=\Phi.
\end{equation}
\end{lemma}

\begin{proof}
By Lemma~\ref{lem:base} we have $\theta J_k'+J_k=J_{k-1}$, that is,
$(I+\theta D)J_k=J_{k-1}$ for all $k\ge 1$ (with the convention
$J_0=\Phi$). Since $\cN=(I+\theta D)^k$, telescoping gives
\[
\cN J_k=(I+\theta D)^k J_k=(I+\theta D)^{k-1}J_{k-1}=\dots=(I+\theta D)^0 J_0=\Phi .
\]
Similarly, from $I_k'=(I_k-I_{k-1})/\theta$ we obtain
$(I-\theta D)I_k=I_{k-1}$, whence
$\wt\cN I_k=(I-\theta D)^k I_k=I_0=\Phi$.

The same follows by substituting \eqref{eq:iter-J} directly
into $\cN J_k=\sum_{n=0}^{k}\binom{k}{n}\theta^n J_k^{(n)}$
(Lemma~\ref{lem:iter} applies up to $n=k$, and the term $m=k$ in
\eqref{eq:iter-J} at $n=k$ is exactly $J_0=\Phi$, so that no $J_{-1}$
is required):
\[
\begin{aligned}
\cN J_k
&=\sum_{m=0}^{k}J_{k-m}\sum_{n=m}^{k}(-1)^{n+m}
  \binom{k}{n}\binom{n}{m}\\
&=\sum_{m=0}^{k}J_{k-m}\binom{k}{m}
  \sum_{j=0}^{k-m}(-1)^j\binom{k-m}{j}
 =\sum_{m=0}^{k}J_{k-m}\binom{k}{m}\delta_{m,k}=\Phi ,
\end{aligned}
\]
where we used the identity
$\binom{k}{n}\binom{n}{m}=\binom{k}{m}\binom{k-m}{n-m}$ and
$(1-1)^{k-m}=\delta_{m,k}$. For $\wt\cN I_k$ the factor $(-1)^n$ of
$\wt\cN$ together with the factor $(-1)^m$ of \eqref{eq:iter-I} gives
the same sign $(-1)^{n+m}$, and the computation repeats verbatim.
\end{proof}

\begin{lemma}\label{lem:commute}
The operators $\cN_{k,\theta}$ and $\wt\cN_{k',\theta'}$ commute for
arbitrary parameters, being polynomials in the operator $D=d/du$ with
constant coefficients.
\end{lemma}
\begin{proof}
We write $\cN_{k,\theta}=\sum_{n=0}^{k}a_nD^n$ and
$\wt\cN_{k',\theta'}=\sum_{m=0}^{k'}b_mD^m$, where
$a_n=\binom{k}{n}\theta^n$ and
$b_m=(-1)^m\binom{k'}{m}(\theta')^m$. For any
$f\in C^{k+k'}$,
\[
\cN_{k,\theta}\wt\cN_{k',\theta'}f
=\sum_{n=0}^{k}\sum_{m=0}^{k'}a_nb_mD^{n+m}f
=\wt\cN_{k',\theta'}\cN_{k,\theta}f .
\]
The required smoothness is provided by assumption~(A2) and the
subsequent bootstrap of regularity, so both orders of applying the
operators are legitimate and coincide.
\end{proof}

\section{Reduction to an ODE for all three models}\label{sec:reduction}

In this section we establish the main reduction of the IDE to an ODE.
Throughout we use the Leibniz rule
\begin{equation}\label{eq:Leibniz}
(u^2\Phi'')^{(n)}=u^2\Phi^{(n+2)}+2nu\Phi^{(n+1)}+n(n-1)\Phi^{(n)},\quad
(u\Phi')^{(n)}=u\Phi^{(n+1)}+n\Phi^{(n)}.
\end{equation}

We begin with the non-life insurance model~(A).

\begin{theorem}\label{thm:A}
Under assumption~(A2) the survival function $\Phi$ satisfies on
$(0,\infty)$ a linear homogeneous ODE of order~$k+2$:
\begin{equation}\label{eq:ODE-A}
g_{k+2}(u)\Phi^{(k+2)}+g_{k+1}(u)\Phi^{(k+1)}+\sum_{n=2}^{k}g_n(u)\Phi^{(n)}+g_1(u)\Phi'=0,
\end{equation}
with the coefficients
\begin{align}
g_{k+2}(u)&=\frac{1}{2}\sigma^2\theta^k u^2,\notag\\
g_{k+1}(u)&=\frac{1}{2}\sigma^2 k\theta^{k-1}u^2+\theta^k(\sigma^2 ku+au+c),\notag\\
g_n(u)&=\frac{1}{2}\sigma^2\tbinom{k}{n-2}\theta^{\,n-2}u^2
       +\tbinom{k}{n-1}\theta^{\,n-1}\bigl(\sigma^2(n-1)u+au+c\bigr) \notag\\
      &\quad+\tbinom{k}{n}\theta^{\,n}\bigl(\frac{1}{2}\sigma^2 n(n-1)+an-\alpha\bigr),\qquad 2\le n\le k,\notag\\
g_1(u)&=au+c+k\theta(a-\alpha).\notag
\end{align}
The term in $\Phi^{(0)}$ is absent, which ensures the reduction of the
order to $k+1$ for $G:=\Phi'$.
\end{theorem}

\begin{proof}
Applying $\cN$ to \eqref{eq:compact-A} and using $\cN J_k=\Phi$
(Lemma~\ref{lem:annihilation}), we obtain
\begin{equation}\label{eq:pre-ODE-A}
\cN\cL\Phi-\alpha\cN\Phi+\alpha\Phi=0.
\end{equation}
We expand $\cN\cL\Phi$ via \eqref{eq:Leibniz}:
\begin{multline}\label{eq:cNL-expanded}
\cN\cL\Phi=\sum_{n=0}^{k}\tbinom{k}{n}\theta^n\Bigl[\frac{1}{2}\sigma^2 u^2\Phi^{(n+2)}+(\sigma^2 nu+au+c)\Phi^{(n+1)}\\
+\bigl(\frac{1}{2}\sigma^2 n(n-1)+an\bigr)\Phi^{(n)}\Bigr].
\end{multline}
The term $-\alpha\cN\Phi+\alpha\Phi$ is written as
$-\alpha\sum_{n=1}^{k}\binom{k}{n}\theta^n\Phi^{(n)}$ (the constant
terms cancel). We collect the coefficients of $\Phi^{(j)}$:

\emph{Leading derivatives.} $j=k+2$: the contribution comes only from
$\Phi^{(n+2)}$ at $n=k$, which gives
$\frac{1}{2}\sigma^2\theta^k u^2=g_{k+2}$.

$j=k+1$: the contributions from $\Phi^{(n+2)}$ at $n=k-1$
($\frac{1}{2}\sigma^2 k\theta^{k-1}u^2$) and from $\Phi^{(n+1)}$ at
$n=k$ ($\theta^k(\sigma^2 ku+au+c)$). The sum is $g_{k+1}$.

\emph{Intermediate derivatives.} For $2\le j\le k$ we collect the
contributions from four sources:
\begin{itemize}[leftmargin=*,itemsep=2pt]
  \item from $\Phi^{(j)}$ at $n=j$:
  $\binom{k}{j}\theta^j\bigl(\frac{1}{2}\sigma^2 j(j-1)+aj\bigr)$,
  \item from $\Phi^{(j)}$ at $n=j-1$ (through $\Phi^{(n+1)}$):
  $\binom{k}{j-1}\theta^{\,j-1}(\sigma^2(j-1)u+au+c)$,
  \item from $\Phi^{(j)}$ at $n=j-2$ (through $\Phi^{(n+2)}$):
  $\binom{k}{j-2}\theta^{\,j-2}\cdot\frac{1}{2}\sigma^2 u^2$,
  \item from $-\alpha\cN\Phi$: $-\alpha\binom{k}{j}\theta^j$.
\end{itemize}
We put these together and obtain $g_j$ from the statement of the
theorem.

\emph{The first derivative.} $j=1$: $\Phi^{(1)}$ arises from
$\Phi^{(n+1)}$ at $n=0$ (the contribution
$\binom{k}{0}\theta^0(au+c)=au+c$), from $\Phi^{(n)}$ at $n=1$ (the
contribution $\binom{k}{1}\theta^1\cdot a=k\theta a$), and from
$-\alpha\cN$ at $n=1$ (the contribution $-\alpha k\theta$). In total
$g_1=au+c+k\theta(a-\alpha)$.

\emph{The term in $\Phi^{(0)}$.} $j=0$: the coefficient of $\Phi$ is
$-\alpha\binom{k}{0}\theta^0+\alpha=0$, so the term in $\Phi$ is absent, and the ODE reduces: $G=\Phi'$
satisfies an equation of order $k+1$ obtained by the substitution
$\Phi^{(j)}=G^{(j-1)}$.
\end{proof}

We now turn to the annuity model~(B): the premium is negative
there, $c<0$, and the jumps are directed upward.

\begin{theorem}\label{thm:B}
Under assumption~(A2), for $c<0$ the function $\Phi$ satisfies on
$(0,\infty)$ a linear homogeneous ODE of order $k+2$:
\begin{equation}\label{eq:ODE-B}
\wt g_{k+2}(u)\Phi^{(k+2)}+\wt g_{k+1}(u)\Phi^{(k+1)}+\sum_{m=2}^{k}\wt g_m(u)\Phi^{(m)}+\wt g_1(u)\Phi'=0,
\end{equation}
with the coefficients (after an overall multiplication by $(-1)^k$ to
normalize the sign of the leading coefficient):
\begin{align}
\wt g_{k+2}(u)&=\frac{1}{2}\sigma^2\theta^k u^2,\notag\\
\wt g_{k+1}(u)&=\theta^k(au+c)+\frac{1}{2} k\theta^{k-1}\sigma^2\,u\,(2\theta-u),\notag\\
\wt g_m(u)&=\frac{1}{2}\sigma^2\theta^{m-2}(-1)^{k-m}\Bigl(
            \tbinom{k}{m-2}u^2-2\tbinom{k}{m-1}\theta(m-1)u\notag\\
          &\qquad\qquad\qquad\qquad\qquad
            +\tbinom{k}{m}\theta^2 m(m-1)\Bigr) \notag\\
          &\quad+(-1)^{k-m+1}\theta^{m-1}a\Bigl(\tbinom{k}{m-1}u-\tbinom{k}{m}\theta m\Bigr) \notag\\
          &\quad+(-1)^{k-m+1}\theta^{m-1}\Bigl(\tbinom{k}{m-1}c+\tbinom{k}{m}\theta\alpha\Bigr),\qquad 2\le m\le k,\notag\\
\wt g_1(u)&=(-1)^k\bigl(au+c+k\theta(\alpha-a)\bigr).\notag
\end{align}
\end{theorem}

\begin{proof}
Applying $\wt\cN$ to \eqref{eq:compact-B} and using
$\wt\cN I_k=\Phi$, we obtain
$\wt\cN\cL\Phi-\alpha\wt\cN\Phi+\alpha\Phi=0$. The expansion differs
from the proof of Theorem~\ref{thm:A} only by the replacement
$\binom{k}{n}\mapsto(-1)^n\binom{k}{n}$ in each term. The coefficient
of $\Phi^{(j)}$ for $2\le j\le k$ takes the form
\begin{multline}\label{eq:Tj}
T_j=(-1)^{j-2}\tbinom{k}{j-2}\theta^{\,j-2}\frac{1}{2}\sigma^2 u^2
    +(-1)^{j-1}\tbinom{k}{j-1}\theta^{\,j-1}\bigl(\sigma^2(j-1)u+au+c\bigr)\\
    +(-1)^j\tbinom{k}{j}\theta^{\,j}\bigl(\frac{1}{2}\sigma^2 j(j-1)+aj-\alpha\bigr).
\end{multline}
Factoring out $(-1)^j$ and multiplying the whole equation by $(-1)^k$
(to normalize the sign of the leading coefficient), we arrive at the
formula for $\wt g_j$ equivalent to the one above. For $j=k+1$ we
similarly obtain
\[
\wt g_{k+1}(u)=\theta^k(au+c)+\frac{1}{2} k\theta^{k-1}\sigma^2 u(2\theta-u),
\]
where the expression $k\theta^{k-1}\sigma^2 u(2\theta-u)/2$ arises by
combining two terms with opposite signs:
$-k\theta^{k-1}\sigma^2 u^2/2+\theta^k\sigma^2 ku=
k\theta^{k-1}\sigma^2 u(2\theta-u)/2$.
\end{proof}

Finally, we consider the mixed model~(C).

\begin{theorem}\label{thm:C}
Under assumption~(A2) the function $\Phi$ satisfies on $(0,\infty)$ a
linear homogeneous ODE of order $k_1+k_2+2$, obtained by applying the
composition of the operators $\wt\cN_{k_1,\theta_1}\,\cN_{k_2,\theta_2}$
to the IDE \eqref{eq:compact-C}. In operator form:
\begin{equation}\label{eq:ODE-C-op}
\wt\cN_{k_1,\theta_1}\cN_{k_2,\theta_2}\cL\Phi+\alpha_1\,\cN_{k_2,\theta_2}\Phi+\alpha_2\,\wt\cN_{k_1,\theta_1}\Phi-(\alpha_1+\alpha_2)\,\wt\cN_{k_1,\theta_1}\cN_{k_2,\theta_2}\Phi=0
\end{equation}
(below we omit the indices of the operators). In the special case
$k_1=k_2=k$, $\theta_1=\theta_2=\theta$ the expanded form reads
\begin{multline}\label{eq:ODE-C}
\sum_{n,m=0}^{k}(-1)^n\tbinom{k}{n}\tbinom{k}{m}\theta^{n+m}\Bigl[\frac{1}{2}\sigma^2 u^2\Phi^{(\ell+2)}+\bigl(\sigma^2 \ell u+au+c\bigr)\Phi^{(\ell+1)}\\
+\bigl(\frac{1}{2}\sigma^2 \ell(\ell-1)+a\ell-(\alpha_1+\alpha_2)\bigr)\Phi^{(\ell)}\Bigr]+\sum_{n=0}^{k}\tbinom{k}{n}\theta^n\bigl(\alpha_1+(-1)^n\alpha_2\bigr)\Phi^{(n)}=0,
\end{multline}
where $\ell:=n+m$. For $k_1=k_2=k$, $\theta_1\neq\theta_2$ it suffices
to replace $\theta^{n+m}$ by $\theta_1^{\,n}\theta_2^{\,m}$, and
$\theta^n(\alpha_1+(-1)^n\alpha_2)$ by
$\alpha_1\theta_2^{\,n}+(-1)^n\alpha_2\theta_1^{\,n}$ (precisely this
form is taken by the coefficients of equation (6) in
\cite{Antipov2026}).
\end{theorem}
\begin{proof}
Applying $\wt\cN\cN$ to \eqref{eq:compact-C}:
\[
\wt\cN\cN\cL\Phi+\alpha_1\wt\cN\cN I_{k_1}+\alpha_2\wt\cN\cN J_{k_2}-(\alpha_1+\alpha_2)\wt\cN\cN\Phi=0.
\]
By Lemmas~\ref{lem:annihilation} and \ref{lem:commute},
\[
\wt\cN\cN I_{k_1}=\cN\wt\cN I_{k_1}=\cN\Phi,\qquad \wt\cN\cN J_{k_2}=\wt\cN\Phi,
\]
whence \eqref{eq:ODE-C-op} follows. In the special case
$k_1=k_2=k$, $\theta_1=\theta_2=\theta$ the composition expands
explicitly:
$\wt\cN\cN f=\sum_{n,m=0}^{k}(-1)^n\binom{k}{n}\binom{k}{m}\theta^{n+m} f^{(n+m)}$.
Applying it to $\cL\Phi$ with the use of \eqref{eq:Leibniz} (with
$\ell=n+m$ in place of $n$) gives the first double sum in
\eqref{eq:ODE-C}, and the replacement
$a\ell\mapsto a\ell-(\alpha_1+\alpha_2)$ inside the bracket accounts
for the term $-(\alpha_1+\alpha_2)\wt\cN\cN\Phi$. Finally,
$\alpha_1\cN\Phi+\alpha_2\wt\cN\Phi$ collapses into the single sum
$\sum_n\binom{k}{n}\theta^n(\alpha_1+(-1)^n\alpha_2)\Phi^{(n)}$.
\end{proof}

In the general case $k_1\neq k_2$ or $\theta_1\neq\theta_2$ the
operators $\cN_{k_2,\theta_2}$ and $\wt\cN_{k_1,\theta_1}$ still commute
as polynomials in $d/du$, and the order of the resulting ODE is
$k_1+k_2+2$. Note the crossed arrangement of the intensities in
\eqref{eq:ODE-C-op}: the intensity $\alpha_1$ (upward jumps, the annihilator
$\wt\cN_{k_1,\theta_1}$) multiplies the operator $\cN_{k_2,\theta_2}$, and
$\alpha_2$ multiplies the operator $\wt\cN_{k_1,\theta_1}$, since
$\wt\cN\cN I_{k_1}=\cN(\wt\cN I_{k_1})=\cN\Phi$. The part of the
equation without $\cL$, i.e.,
$-(\alpha_1+\alpha_2)\wt\cN\cN\Phi+\alpha_1\cN\Phi+\alpha_2\wt\cN\Phi$,
does not contain $\Phi^{(0)}$ (the constant terms of all three
operators are equal to $1$). Hence, as in the one-sided models,
equation \eqref{eq:ODE-C-op} is an equation of order $k_1+k_2+1$ for
$G=\Phi'$.

\begin{remark}\label{rem:selection}
Theorems~\ref{thm:A}--\ref{thm:C} establish that the survival function
$\Phi$ satisfies the corresponding ODE. The converse inclusion does not
follow automatically: the solution space of a linear ODE of order
$k+2$ is wider than the original probabilistic problem. The physical
solution additionally satisfies conditions consistent with the nature
of $\Phi$:
\begin{enumerate}[label=(\roman*),leftmargin=2em]
  \item $\Phi(z)\equiv 0$ for $z<0$,
  \item $\lim_{u\to\infty}\Phi(u)=1$ (from the non-degeneracy
  condition~(A1)),
  \item $\Phi(0+)\in(0,1)$ in model~(A) and in model~(C) for $c\ge0$
  (the positivity follows from Proposition~\ref{prop:C-zero} and
  Remark~\ref{rem:c0}, and $\Phi(0+)<1$ because $\Psi(0+)>0$),
  and $\Phi(0+)=0$ in model~(B) and in model~(C) for $c<0$ (property
  (iii) of Sect.~\ref{sec:model}),
  \item $G=\Phi'\ge 0$ (the monotonicity of $\Phi$),
  \item $G\in L^1((0,\infty))$ (the integrability of the density,
  equivalent to $1-\Phi(0+)=\int_0^\infty G\,du<\infty$).
\end{enumerate}
Conditions (iv)--(v) narrow the class of admissible solutions but do
not by themselves prove uniqueness: in model~(A) they do not separate
the power branch from the exponentially small modes (all of them are
integrable). The following proposition gives an exact criterion for an ODE
solution to satisfy the original IDE. For $k=1$ it
reduces to the condition $c\Phi'(0+)=\alpha\Phi(0+)$ which closes the
verification theorem in \cite{PRYu2026}.
\end{remark}

\begin{proposition}[converse criterion]\label{prop:converse}
Let $\Phi\in C^{K+2}((0,\infty))$ be bounded, where $K$ is the order
of the annihilator ($K=k$ in models~(A), (B), and $K=k_1+k_2$ in
model~(C)), and let $\cR$ be the residual of the IDE:
\begin{align*}
\cR&:=\cL\Phi+\alpha J_k-\alpha\Phi &&\text{in model~(A)},\\
\cR&:=\cL\Phi+\alpha I_k-\alpha\Phi &&\text{in model~(B)},\\
\cR&:=\cL\Phi+\alpha_1I_{k_1}+\alpha_2J_{k_2}-(\alpha_1+\alpha_2)\Phi &&\text{in model~(C)}.
\end{align*}
Then $\Phi$ satisfies the reduced ODE
(Theorems~\ref{thm:A}--\ref{thm:C}) if and only if $\cR$ lies in the
kernel of the annihilator $\cA$ ($\cA=\cN$, $\wt\cN$, $\wt\cN\cN$
respectively), i.e., is a quasi-polynomial: in~(A)
$\cR=e^{-u/\theta}Q(u)$ with $\deg Q\le k-1$, in~(B)
$\cR=e^{u/\theta}Q(u)$ with $\deg Q\le k-1$, and in~(C)
$\cR=e^{u/\theta_1}Q_1(u)+e^{-u/\theta_2}Q_2(u)$ with
$\deg Q_i\le k_i-1$. In particular, a solution of the reduced ODE
satisfies the IDE if and only if
\begin{equation}\label{eq:converse}
\lim_{u\to0+}\cR^{(j)}(u)=0,\qquad j=0,1,\dots,K-1,
\end{equation}
and these limits exist for every bounded solution of class
$C^{K+2}((0,\infty))$ of the ODE.
\end{proposition}

\begin{proof}
Lemmas~\ref{lem:annihilation} and \ref{lem:commute} are identities valid
for any bounded continuous $\Phi$, so $\cA \cR$ coincides with the
left-hand side of the reduced ODE: for example, in model~(A)
$\cN \cR=\cN\cL\Phi+\alpha\Phi-\alpha\cN\Phi$, and in model~(C)
$\wt\cN\cN \cR=\wt\cN\cN\cL\Phi+\alpha_1\cN\Phi+\alpha_2\wt\cN\Phi-(\alpha_1+\alpha_2)\wt\cN\cN\Phi$.
The ODE is therefore equivalent to $\cA \cR=0$, i.e.,
$\cR\in\ker\cA$, and the kernel of an operator of order $K$ with
constant coefficients consists precisely of the indicated
quasi-polynomials. Every element of the kernel is an entire function
uniquely determined by the values $\cR^{(j)}(0)$, $j<K$ (the
uniqueness of the solution of the Cauchy problem for $\cA \cR=0$).
Hence $\cR\equiv0$ if and only if \eqref{eq:converse} holds. The
existence of the limits follows from the fact that on $(0,\infty)$
the function $\cR$ coincides with an entire function.
\end{proof}

\begin{remark}\label{rem:converse}
Conditions \eqref{eq:converse} are expressed through $\Phi$ and the
convolutions. In model~(A) we have $J_k(u)=O(u^k)$ and, by
Lemma~\ref{lem:iter}, $J_k^{(j)}(0+)=0$ for $j<k$, so that
\eqref{eq:converse} takes the form
$\lim_{u\to0+}[(\cL-\alpha)\Phi]^{(j)}(u)=0$, $j<k$. For $k=1$ and
$u^2\Phi''\to0$ this is the relation $c\Phi'(0+)=\alpha\Phi(0+)$,
which coincides with the condition of equivalence of the ODE and the
IDE in \cite{PRYu2026}. In models~(B), (C) the conditions involve the
values $I_k^{(j)}(0+)$, i.e., functionals of $\Phi$ on the whole
half-line. The number of effective conditions, however, does not
exceed $K$ and depends on the model. For a bounded solution of the
ODE the residual $\cR$ is bounded as $u\to\infty$: by
Lemma~\ref{lem:actual} and the proof of
Proposition~\ref{prop:realization} in Sect.~\ref{sec:asymptotics}
(the exclusion of the growing modes), such a solution is a combination
of a constant, the power branch, and exponentially decaying modes with
termwise differentiable asymptotics, so that $\cL\Phi\to0$ and the
convolutions are bounded. The growing quasi-polynomials
$e^{u/\theta}Q(u)$ in~(B) and $e^{u/\theta_1}Q_1(u)$ in~(C) cannot be
bounded. In model~(B), then, every bounded solution of the ODE
satisfies the IDE (all $k$ conditions hold automatically). In
model~(C), conditions \eqref{eq:converse} reduce to at most $k_2$
conditions $Q_2\equiv0$ (the coefficients of the polynomial $Q_2$),
and in model~(A) to at most $k$ conditions $Q\equiv0$. The
independence of these conditions is not used below. In model~(A) all
$k$ conditions are independent: solutions bounded at zero have
asymptotic series $\sum_{n\ge0}b_nu^n$ there with free
$b_0,\dots,b_k$, because the equation of Theorem~\ref{thm:A} turns out to have at zero the
   same structure as equation \eqref{eq:ODE-C-op} in
Proposition~\ref{prop:C-zero} (the leading coefficient
$\frac{1}{2}\sigma^2\theta^ku^2$, the coefficient of $\Phi^{(k+1)}$
equal to $c\theta^k\neq0$ at zero, the indicial equation
$\nu(\nu-1)\cdots(\nu-k)=0$, and the mode of the form
\eqref{eq:C-irr} excluded for $c>0$ by boundedness), so that the proof
together with the asymptotic existence theorem at zero carries over
verbatim. Here
$[(\cL-\alpha)\Phi]^{(j)}(0+)=j!\,[c(j+1)b_{j+1}+(\frac{1}{2}\sigma^2j(j-1)+aj-\alpha)b_j]$,
so that the conditions form a triangular system in $b_1,\dots,b_k$
with a nonzero diagonal $c\,(j+1)!$. In model~(C) the conditions are
independent in the numerically tested examples. The dimension $K$ of
the kernel of the annihilator agrees with the structure of the kernel
in \cite[Sect.~7]{KPuh2020}: for $k_1=k_2=1$ it is shown there that
the kernel of the ODE operator is four-dimensional, with a basis
consisting of the constant, the function $\Phi$, and two solutions
that do not satisfy the IDE, while the kernel of the annihilator (the
operator $T$ in the notation of \cite{KPuh2020}) is two-dimensional,
in agreement with $K=2$.
\end{remark}

\section{Asymptotic analysis}\label{sec:asymptotics}

\subsection{The case $k=2$: formal power branch}\label{subsec:k2}

For $k=2$ the ODE for $G:=\Phi'$ is of the third order:
\begin{equation}\label{eq:ODE-k2}
\frac{1}{2}\sigma^2\theta^2 u^2 G'''+g_3(u)G''+g_2(u)G'+g_1(u)G=0,
\end{equation}
with
\begin{align*}
g_3(u)&=\sigma^2\theta u^2+\theta^2(2\sigma^2 u+au+c),\\
g_2(u)&=\frac{1}{2}\sigma^2 u^2+2\theta(\sigma^2 u+au+c)+\theta^2(\sigma^2+2a-\alpha),\\
g_1(u)&=2\theta a+au+c-2\theta\alpha.
\end{align*}

In normalized form ($G'''+q_3 G''+q_2 G'+q_1 G=0$, where
$g_4:=g_{k+2}|_{k=2}=\frac{1}{2}\sigma^2\theta^2u^2$), the expansions
of the coefficients in powers of $1/u$ have the explicit form:
\begin{align}
q_3(u)&=\frac{g_3(u)}{g_4(u)}
       =\frac{2}{\theta}+\frac{2(a+2\sigma^2)}{\sigma^2}\cdot\frac1u+\frac{2c}{\sigma^2}\cdot\frac1{u^2},\\
q_2(u)&=\frac{g_2(u)}{g_4(u)}
       =\frac1{\theta^2}+\frac{4(\sigma^2+a)}{\sigma^2\theta}\cdot\frac1u
        +\frac{2\bigl(2c+\theta(\sigma^2+2a-\alpha)\bigr)}{\theta\sigma^2}\cdot\frac1{u^2},\\
q_1(u)&=\frac{g_1(u)}{g_4(u)}
       =\frac{2a}{\sigma^2\theta^2}\cdot\frac1u+\frac{2(2\theta a+c-2\theta\alpha)}{\sigma^2\theta^2}\cdot\frac1{u^2}.
\end{align}
For brevity, we write
\[
q_3(u)=\frac{2}{\theta}+c_{31}\frac1u+c_{32}\frac1{u^2},\quad
q_2(u)=\frac1{\theta^2}+c_{21}\frac1u+c_{22}\frac1{u^2},\quad
q_1(u)=c_{11}\frac1u+c_{12}\frac1{u^2}.
\]

\begin{theorem}[formal power branch, $k=2$]\label{thm:formal-k2}
Let $k=2$ (model~(A)). Equation \eqref{eq:ODE-k2} admits a unique
formal solution of power type (a Thom\'e series)
\begin{equation}\label{eq:Gformal}
G_{\mathrm f}(u)=u^{-p}\sum_{j\ge 0}a_j u^{-j},\qquad a_0=1,\qquad p=\beta+1=\frac{2a}{\sigma^2}>1,
\end{equation}
where the exponent $p$ is determined by the dominant balance, and the
coefficients $(a_j)_{j\ge 1}$ are uniquely determined by the recurrence
\begin{multline}\label{eq:rec}
0=a_j\Bigl(c_{11}-\frac{p+j}{\theta^2}\Bigr)\\
  +a_{j-1}\Bigl(c_{12}-c_{21}(p{+}j{-}1)+\frac{2}{\theta}(p{+}j{-}1)(p{+}j{-}2)+\frac{4}{\theta}(p{+}j{-}1)\Bigr)\\
  +a_{j-2}\Bigl(-c_{22}(p{+}j{-}2)+2c_{31}(p{+}j{-}2)-6(p{+}j{-}2)-6(p{+}j{-}2)(p{+}j{-}3)\\
  \hspace{3em}+c_{31}(p{+}j{-}2)(p{+}j{-}3)-(p{+}j{-}2)(p{+}j{-}3)(p{+}j{-}4)\Bigr)\\
  +a_{j-3}\Bigl(2c_{32}(p{+}j{-}3)+c_{32}(p{+}j{-}3)(p{+}j{-}4)\Bigr),
\end{multline}
with the notation $c_{ij}$ from the expansions of $q_i(u)$. The
convention $a_n=0$ for $n<0$ means that for $j=1,2$ some terms of
\eqref{eq:rec} are absent. For clarity we display the first two steps
separately:
\begin{align}\label{eq:rec-j1}
0&=a_1\Bigl(c_{11}-\frac{p+1}{\theta^2}\Bigr)
 +\Bigl(c_{12}-c_{21}p+\frac{2}{\theta}p(p-1)+\frac{4}{\theta}p\Bigr),\\
\label{eq:rec-j2}
0&=a_2\Bigl(c_{11}-\frac{p+2}{\theta^2}\Bigr)
 +a_1\Bigl(c_{12}-c_{21}(p+1)+\frac{2}{\theta}(p+1)p+\frac{4}{\theta}(p+1)\Bigr)\notag\\
&\quad+\Bigl(-c_{22}p+2c_{31}p-6p-6p(p-1)+c_{31}p(p-1)-p(p-1)(p-2)\Bigr),
\end{align}
while the general formula \eqref{eq:rec} applies for $j\ge3$. The step
$j=0$ of \eqref{eq:rec} is the indicial equality $c_{11}=p/\theta^2$.
Hence the factor multiplying $a_j$ equals
$c_{11}-(p+j)/\theta^2=-j/\theta^2\ne0$ for all $j\ge1$, and the
recursion is solvable uniquely. The series \eqref{eq:Gformal} is, in
general, divergent (a Gevrey class of order $1$), but it serves as the
asymptotic expansion of the actual solution $G_{\mathrm{pow}}$ of
equation \eqref{eq:ODE-k2} (Lemma~\ref{lem:actual}), and the physical
density $G=\Phi'$ contains this branch with a coefficient $\varpi>0$
(Proposition~\ref{prop:realization}), so that
\begin{equation}\label{eq:asymp}
\Psi(u)\sim C u^{-\beta},\qquad u\to\infty,\qquad C=\frac{\varpi}{\beta}>0.
\end{equation}
\end{theorem}

\begin{proof}
\emph{Step 1. The nature of the singular point.} We substitute
$z=1/u$, $w(z):=G(1/z)$, which converts \eqref{eq:ODE-k2} to the equation
\begin{equation}\label{eq:z-form}
-z^6 w'''+(-6z^5+\wt q_3 z^4)w''+(-6z^4+2\wt q_3 z^3-\wt q_2 z^2)w'+\wt q_1 w=0,
\end{equation}
where $\wt q_i(z):=q_i(1/z)$. Since
$\wt q_3(z)=2/\theta+O(z)$, $\wt q_2(z)=1/\theta^2+O(z)$, and
$\wt q_1(z)=O(z)$ as $z\to0$, the point $z=0$ (corresponding to
$u=\infty$) is an irregular singular point of rank~1: the
characteristic equation for the exponential solutions $e^{\mu u}$ of
\eqref{eq:ODE-k2} has the form $\mu(\mu+1/\theta)^2=0$, and the root
$\mu=0$, corresponding to the power branch, is simple (see
\cite[Sects.~12--15]{Wasow1965}).

\emph{Step 2. The indicial equation.} We substitute the power expansion
\[
G(u)=u^{-p}\sum_{j\ge 0}a_j u^{-j},\qquad a_0=1,
\]
into the normalized form of \eqref{eq:ODE-k2}. As $u\to\infty$ we
obtain the following orders of the contributions:
\[
\begin{array}{c|c}
\text{Term of the equation} & \text{Leading contribution as }u\to\infty \\ \hline
G''' & -a_0\,p(p+1)(p+2)\,u^{-p-3} \\[2pt]
q_3 G'' & \displaystyle\frac{2}{\theta}\cdot a_0\,p(p+1)\,u^{-p-2} \\[2pt]
q_2 G'   & \displaystyle\frac{1}{\theta^2}\cdot(-a_0 p)\,u^{-p-1} \\[2pt]
q_1 G    & \displaystyle c_{11}\cdot a_0\,u^{-p-1}
\end{array}
\]

The dominant balance (the highest power of $u$) occurs at
$u^{-p-1}$: comparing $u^{-p-1}$ with $u^{-p-2}$ and $u^{-p-3}$ as
$u\to\infty$, we find that the leading terms are $q_2 G'$ (the constant part
$1/\theta^2$) and $q_1 G$ (the term $c_{11}/u$). All other
contributions ($G'''$, the subdominant parts of $q_3 G''$ and
$q_2 G'$) belong to higher powers of $1/u$ and play no role in
determining $p$.

We equate to zero the coefficient $C_{p+1}$ at $u^{-p-1}$ (with
$a_0=1$) and obtain
\[
C_{p+1}=c_{11}-\frac{p}{\theta^2}=0\quad\Longrightarrow\quad p=c_{11}\theta^2=\frac{2a}{\sigma^2}=\beta+1.
\]

\emph{Step 3. The recurrence.} We systematically count all contributions
to the coefficient of $u^{-p-j-1}$ from each of the four terms of
\eqref{eq:ODE-k2} in normalized form, which gives \eqref{eq:rec}. With
$P:=p+j$ the contributions of the four terms are
\begin{align*}
G'''&:\ -a_{j-2}(P-2)(P-1)P,\\
q_3G''&:\ \frac2\theta a_{j-1}(P-1)P+c_{31}a_{j-2}(P-2)(P-1)+c_{32}a_{j-3}(P-3)(P-2),\\
q_2G'&:\ -\theta^{-2}a_jP-c_{21}a_{j-1}(P-1)-c_{22}a_{j-2}(P-2),\\
q_1G&:\ c_{11}a_j+c_{12}a_{j-1}.
\end{align*}
Expanding the products, we arrive at the form \eqref{eq:rec}.

\emph{Step 4. The asymptotics of $\Psi$.} By Lemma~\ref{lem:actual},
to the formal series \eqref{eq:Gformal} there corresponds an actual
solution $G_{\mathrm{pow}}$, and the physical density has the form
$G=\varpi\,G_{\mathrm{pow}}+O(e^{-\delta u})$ for some $\delta>0$ (the
remaining solutions are exponentially small). By
Proposition~\ref{prop:realization}, $\varpi\ne0$, and from $G\ge0$ (the
monotonicity of $\Phi$) we get $\varpi>0$. From $\Psi(\infty)=0$, by
integration,
\[
\Psi(u)=\int_u^\infty G(v)\,dv\sim \frac{\varpi}{p-1}u^{1-p}
=\frac{\varpi}{\beta}u^{-\beta}.
\]
We fix the normalization $a_0=1$ only for the formal series: the
amplitude $\varpi$ (the connection coefficient) is determined by the
original problem, not by the reduced equation.
\end{proof}

\begin{remark}\label{rem:dominant}
The exponent $p$ should not be extracted from the coefficient of the
lowest power $u^{-p-3}$, which is occupied by the term $G'''$: at this power the
coefficients $a_1,a_2$ of the series appear along with $a_0$, so the
corresponding relation enters the recurrence and does not constrain
$p$. For an irregular singular point, the equation for the exponent of
the power branch is given by the dominant balance: since
$-p-1>-p-2>-p-3$, it is formed by the terms $q_2G'$ and $q_1G$, while
the highest derivative $G'''$ is subdominant. The form \eqref{eq:rec}
is equivalent to the compact one: with $P:=p+j$ the factors multiplying
$a_{j-1}$, $a_{j-2}$, $a_{j-3}$ equal
$c_{12}-c_{21}(P-1)+\frac2\theta(P-1)P$,
$-c_{22}(P-2)+c_{31}(P-2)(P-1)-(P-2)(P-1)P$ and $c_{32}(P-3)(P-2)$
respectively. For example, with $a_0=1$ we get
$a_1=2(c-2\alpha\theta)/\sigma^2$ and
$a_2=-2(3a\alpha\theta^2-4\alpha^2\theta^2+4\alpha c\theta-c^2)/\sigma^4$.
\end{remark}

\subsection{The exponent of the power branch for a general $k$}\label{subsec:generalk}
In normalized form the ODE for $G=\Phi'$ for general $k$ reads
\[G^{(k+1)}+q_{k+1}(u)G^{(k)}+\dots+q_2(u)G'+q_1(u)G=0.\]
The coefficients $q_1$, $q_2$ have the expansions
\begin{align}
q_1(u)&=\frac{2a}{\sigma^2\theta^k}\cdot\frac1u+\frac{2(\theta ka+c-\theta k\alpha)}{\sigma^2\theta^k}\cdot\frac1{u^2},\\
q_2(u)&=\frac1{\theta^k}+\frac{2k(\sigma^2+a)}{\sigma^2\theta^{k-1}}\cdot\frac1u
       +\frac{k(k-1)\theta(\sigma^2+2a-\alpha)+2kc}{\sigma^2\theta^{k-1}}\cdot\frac1{u^2},
\end{align}
obtained by the direct division $g_n/g_{k+2}$ (taking into account
$g_{k+2}=\frac{1}{2}\sigma^2\theta^k u^2$). Here and below
$c_{11}:=2a/(\sigma^2\theta^k)$ is the coefficient of $1/u$ in $q_1$
(for $k=2$ this is the previous $c_{11}$). The formal power expansion
(a Thom\'e series) $G(u)=u^{-p}\sum_{j\ge 0}a_j u^{-j}$ with $a_0=1$
gives $G^{(i)}(u)\sim(-1)^i\,p(p+1)\cdots(p+i-1)\,u^{-p-i}$ as
$u\to\infty$. The dominant balance now reads as follows:

\begin{center}
\small\setlength{\tabcolsep}{4pt}
\begin{tabular}{c|c|l}
Term & Leading term as $u\to\infty$ & Order \\ \hline
$G^{(k+1)}$    & $a_0\,(-1)^{k+1}p(p+1)\cdots(p+k)\,u^{-p-(k+1)}$ & $u^{-p-(k+1)}$ \\[2pt]
$q_{m} G^{(m-1)}$, $3\le m\le k+1$ & $\displaystyle q_{m0}\cdot a_0\,(-1)^{m-1}p\cdots(p+m-2)\,u^{-p-(m-1)}$ & $u^{-p-(m-1)}$ \\[2pt]
$q_2 G'$        & $\displaystyle\frac{1}{\theta^k}\cdot(-a_0 p)\,u^{-p-1}$ & $\boldsymbol{u^{-p-1}}$ \\[2pt]
$q_1 G$         & $\displaystyle c_{11}\cdot a_0\,u^{-p-1}$ & $\boldsymbol{u^{-p-1}}$
\end{tabular}
\end{center}

Here $q_{m0}$ are the leading constant values of
$q_m(u)=g_m/g_{k+2}$ as $u\to\infty$. For $2\le m\le k+1$ they are
finite, because the $g_m$ are polynomials of degree at most $2$, while
$g_{k+2}\propto u^2$.

The highest power of $u$ among all terms is $u^{-p-1}$ (since
$-p-1 > -p-m+1$ for $m\ge 3$). The dominant balance is therefore
determined exclusively by $q_2 G'$ and $q_1 G$:
\[
\frac{1}{\theta^k}\cdot(-p) + c_{11} = 0 \quad\Longrightarrow\quad p = c_{11}\theta^k = \frac{2a}{\sigma^2} = \beta+1.
\]

All the remaining terms ($G^{(k+1)}$ and $q_m G^{(m-1)}$ for $m\ge 3$)
contribute $O(u^{-p-2})$ and do not enter the determining equation.
The exponent $p$ thus does not depend on $k$, $c$, $\alpha$, $\theta$:
it is determined solely by the parameters of the risky asset
$(a,\sigma)$. For the same reason the term $a_j$ enters the
coefficient of $u^{-p-j-1}$ only through $q_2G'$ and $q_1G$, so the
factor multiplying $a_j$ in the recurrence equals
$c_{11}-(p+j)/\theta^k=-j/\theta^k\neq0$, and the formal series is
determined uniquely for any $k$. The step $j=1$ of the same recurrence
(the coefficient of $u^{-p-2}$ is assembled from $q_3G''$, $q_2G'$ and
$q_1G$) gives
$a_1=\theta^k\bigl[q_{30}\,p(p+1)-q_{21}\,p+c_{12}\bigr]$, where
$q_{30}=k/\theta^{k-1}$ is the limit of $q_3=g_3/g_{k+2}$ as
$u\to\infty$ (for $k=1$ the role of $q_3$ is played by the unit
coefficient of $G''$), while
$q_{21}=2k(\sigma^2+a)/(\sigma^2\theta^{k-1})$ and
$c_{12}=2(\theta ka+c-\theta k\alpha)/(\sigma^2\theta^k)$ are the
coefficients of $1/u$ in $q_2$ and of $1/u^2$ in $q_1$. Substituting
$p=2a/\sigma^2$, we get $a_1=2(c-k\alpha\theta)/\sigma^2$. We have thus proved the
following proposition.

\begin{proposition}[formal power branch for a general $k$]\label{prop:general-k}
For model~(A) with an arbitrary $k\ge 1$ the reduced ODE for $G=\Phi'$
admits a unique formal solution of the form
$u^{-p}\sum_{j\ge0}a_ju^{-j}$ with $a_0=1$ and the exponent
$p=\beta+1=2a/\sigma^2$, which does not depend on $k$, $c$, $\alpha$
or $\theta$. The factor multiplying $a_j$ in the recurrence is
$-j/\theta^k\neq0$, and the first coefficient is
$a_1=2(c-k\alpha\theta)/\sigma^2$. Proposition~\ref{prop:realization} establishes
the realization of this branch in the physical solution and the
asymptotics $\Psi(u)\sim C_ku^{-\beta}$, $C_k>0$.
\end{proposition}

The irregular singular point $u=\infty$ of rank~1 generates a
fundamental set of solutions consisting of the power solution and
exponentially small solutions (formal Thom\'e series, see
\cite{Wasow1965}, and for second-order equations
\cite[Chapter~7]{Olver1974}, \cite[Sect.~3.1]{Slavyanov2000}): the
characteristic polynomial of the equation for $G$ is
$\mu(1+\theta\mu)^k$, the root $\mu=0$ is simple for all $k$, and the
root $\mu=-1/\theta$ has multiplicity $k$. For example, for $k=1$
(model~(A)) the ODE \eqref{eq:ODE-k1} has two fundamental solutions
near $u=\infty$: the power one $\sim u^{-p}$ and the exponentially
small one $\sim e^{-u/\theta}u^{-2}$. Substituting
$G=e^{\mu u}u^{\rho}$ directly into \eqref{eq:ODE-k1}, we confirm this: the order
$u^2$ gives $\frac{1}{2}\sigma^2\mu(\theta\mu+1)=0$, i.e.,
$\mu=-1/\theta$, and the order $u$ gives
$-\frac{1}{2}\sigma^2(\rho+2)=0$, i.e., $\rho=-2$. For $k=1$, the
connection coefficient of the power branch in \cite{PRYu2026} is
identified with the connection coefficient of the doubly confluent
Heun equation. Its positivity is proved there from a first integral,
and for an integer $2a/\sigma^2$ it is written explicitly in closed
form.

\begin{remark}\label{rem:realization}
The statement $\Psi(u)\sim Cu^{-\beta}$ with $C>0$ is equivalent to a
nonzero coefficient of the power branch in the expansion of the
physical density $G=\Phi'$ in the actual fundamental set
(Lemma~\ref{lem:actual}). For $k=1$ this is established in
\cite{FKP2002} and, from a first integral of the ODE, in
\cite{PRYu2026}. For general claim distributions with
$\bbe[\xi^{\beta}]<\infty$ the asymptotics with $C>0$ were obtained by
the integral-equation method in \cite[Theorem~1]{P05143}. Here, for
all $k$, the nonvanishing follows from the power-law lower bound for
$\Psi$ (Proposition~\ref{prop:realization}). Conditions (iv)--(v) of
Remark~\ref{rem:selection} (the monotonicity and integrability of $G$)
narrow the solution space but do not separate the power branch
$u^{-p}$ from the exponentially small modes (both are integrable).
This is exactly why a separation requires an external argument (the
lower bound). The criterion of Proposition~\ref{prop:converse} selects
solutions of the IDE among solutions of the ODE, but by itself it does
not imply the nonvanishing of the coefficient $c_\beta$ at the power
branch (Proposition~\ref{prop:realization}). Borel summability of the
series \eqref{eq:Gformal} is not needed for the asymptotics. For
$k=2$, in addition, the coefficients for generic parameters (with nonzero Stokes
multipliers relating $G_{\mathrm{pow}}$ to the exponentially small
modes) grow factorially,
$\limsup_{j\to\infty}(|a_j|/j!)^{1/j}=\theta$ (for
$8\alpha\le\sigma^2$ also $a_j/(j\,a_{j-1})\to\theta$, whereas for
$8\alpha>\sigma^2$ the exponents $s$ of the exponentially small modes
$e^{-u/\theta}u^{s}$, $\sigma^2(s+1)(s+2)+2\alpha=0$, are complex, and
$a_j$ changes sign infinitely often --- log-periodically, i.e., more
and more rarely --- and all this is confirmed numerically), so that
the Borel transform of the series has a singularity at the point
$1/\theta$ of the positive half-axis, the direction $u\to+\infty$ is a
Stokes line, and the formal series determines the actual solution only
up to exponentially small terms $e^{-u/\theta+o(u)}$.
\end{remark}

\subsection{All three models: the exponent, the fundamental set, the tail realization}\label{subsec:allmodels}

The detailed asymptotic analysis above was carried out for model~(A),
for which the most complete basis for comparison is available
\cite{FKP2002,KPuh2020,PRYu2026}. The exponent of the power
asymptotics, however, is the same for all three models.

\begin{proposition}\label{prop:exponent}
In each of the three models the equation for the exponent $s$ of
the algebraic (power) branch $\Phi\sim u^{s}$ of the reduced ODE at the
irregular singular point $u=\infty$ is one and the same:
\begin{equation}\label{eq:indicial-common}
s\Bigl(\frac{1}{2}\sigma^2(s-1)+a\Bigr)=0,
\end{equation}
with roots $s=0$ and $s=1-2a/\sigma^2=-\beta$. The exponent $\beta=2a/\sigma^2-1$ depends neither on the shape
parameter $k$ nor on the parameters of the risk process $c$, $\alpha$,
$\theta$.
\end{proposition}

\begin{proof}
Substitute $\Phi\sim u^{s}$ into the reduced equation. The
annihilating operators take $u^{s}$ again into $u^{s}$ with leading
coefficient $1$, so the coefficient of the highest power $u^{s}$ is
provided by the operator $\cL$ alone and equals
$\frac{1}{2}\sigma^2 s(s-1)+as$. The integral terms contribute to no
such power: the differences $-\alpha\cN\Phi+\alpha\Phi$ in model~(A),
$-\alpha\wt\cN\Phi+\alpha\Phi$ in model~(B), and
$-(\alpha_1+\alpha_2)\wt\cN\cN\Phi+\alpha_1\cN\Phi+\alpha_2\wt\cN\Phi$
in model~(C) start at order $u^{s-1}$. We equate to zero the
coefficient of $u^{s}$ and obtain \eqref{eq:indicial-common}. The argument
does not depend on $k$: all three operators $\cN$, $\wt\cN$,
$\wt\cN\cN$ are polynomials in $D=d/du$ with constant term $1$, so for
$\Phi\sim u^s$ we have $\cN\Phi=\Phi+O(u^{s-1})$ and so on, and the
integral part at order $u^s$ cancels exactly ($-\alpha\cN\Phi+\alpha\Phi=O(u^{s-1})$
in~(A),~(B), and
$-(\alpha_1+\alpha_2)\wt\cN\cN\Phi+\alpha_1\cN\Phi+\alpha_2\wt\cN\Phi=(-(\alpha_1+\alpha_2)+\alpha_1+\alpha_2)\Phi+O(u^{s-1})=O(u^{s-1})$
in~(C)).
\end{proof}

In the annuity model~(B) all $k$ exponential modes grow (the symbol of
the annihilator $(1-\theta\mu)^k$ has the unique root $\mu=1/\theta>0$).
Hence the bounded solutions of the reduced ODE form a two-dimensional
space spanned by $1$ and the power branch, and we obtain the asymptotics
$\Psi(u)\sim Cu^{-\beta}$ from the ODE alone
(Proposition~\ref{prop:realization}), the coefficients of the power
series being given by a recurrence analogous to \eqref{eq:rec}. In the
mixed model~(C) the solution space of the reduced ODE
\eqref{eq:ODE-C-op} of order $k_1+k_2+2$ at the irregular point
$u=\infty$ is wider than in the one-sided models. We introduce the
symbol of the annihilator
\begin{equation}\label{eq:symbolC}
A(\mu):=(1-\theta_1\mu)^{k_1}(1+\theta_2\mu)^{k_2},\qquad\text{so that}\quad \wt\cN_{k_1}\cN_{k_2}\,e^{\mu u}=A(\mu)\,e^{\mu u}.
\end{equation}

\begin{lemma}[fundamental set of model~(C)]\label{lem:fundC}
The $(k_1+k_2+2)$-dimensional space of formal solutions (Thom\'e
series, generalized at repeated roots by fractional powers of
$u$ in the exponent and logarithmic terms) of equation \eqref{eq:ODE-C-op} at
the irregular singular point $u=\infty$ splits into three groups:
\begin{enumerate}[label=\textup{(\alph*)},leftmargin=2.2em,itemsep=1pt]
  \item two algebraic branches $u^{s}$, $s\in\{0,-\beta\}$, generated
  by the diffusion part $\cL$ (Proposition~\ref{prop:exponent}),
  \item $k_1$ solutions with leading exponential factor $e^{u/\theta_1}$
  (exponentially growing), corresponding to the root $\mu=1/\theta_1$
  of the symbol \eqref{eq:symbolC} of multiplicity $k_1$ (upward
  jumps),
  \item $k_2$ solutions with leading exponential factor $e^{-u/\theta_2}$
  (exponentially decaying), corresponding to the root $\mu=-1/\theta_2$
  of multiplicity $k_2$ (downward jumps).
\end{enumerate}
\end{lemma}

\begin{proof}
Here $k_1,k_2\in\bbn$, so that $\wt\cN_{k_1}\cN_{k_2}$ is a
differential operator of order $k_1+k_2$ with constant coefficients.
We substitute $\Phi=e^{\mu u}v$ with $v\sim u^{\rho}$. By the Leibniz
formula,
\[
\wt\cN_{k_1}\cN_{k_2}[e^{\mu u}v]
=e^{\mu u}A(\mu+D)v
=e^{\mu u}\bigl(A(\mu)v+A'(\mu)v'+\cdots\bigr),
\]
while
\[
\cL\Phi=e^{\mu u}\left[\frac{1}{2}\sigma^2\mu^2u^{\rho+2}
+\mu(\sigma^2\rho+a)u^{\rho+1}+O(u^{\rho})\right].
\]
The leading term is
\[
\wt\cN_{k_1}\cN_{k_2}\cL\Phi
=\frac{1}{2}\sigma^2\mu^2A(\mu)u^{\rho+2}e^{\mu u}(1+O(u^{-1})),
\]
whereas the remaining terms of \eqref{eq:ODE-C-op} are of order
$u^{\rho}e^{\mu u}$. For $\mu\neq0$ the condition that the leading term
vanish is equivalent to $A(\mu)=0$, i.e., $\mu=1/\theta_1$ or
$\mu=-1/\theta_2$ with multiplicities $k_1$ and $k_2$. The first family
of modes grows, and the second one decays. To the sector $\mu=0$ there
correspond the two algebraic branches with the equation
$s(\frac{1}{2}\sigma^2(s-1)+a)=0$ of
Proposition~\ref{prop:exponent}.

The algebraic prefactor only refines the subdominant structure. For a
simple root ($k_i=1$: then $A(\mu)=0$ but $A'(\mu)\neq0$ --- for
example, $A'(1/\theta_1)=-\theta_1(1+\theta_2/\theta_1)^{k_2}\neq0$),
the next order $u^{\rho+1}$ gives
$\frac{1}{2}\sigma^2\mu^2(\rho+2)A'(\mu)=0$, whence $\rho=-2$. This
agrees with the branch $e^{-u/\theta}u^{-2}$ of \eqref{eq:ODE-k1} (the
case $k=1$, cf.\ \cite{PRYu2026}), and for $k_1=k_2=1$ we obtain the whole
fundamental set
$\{u^0,\,u^{-\beta},\,e^{u/\theta_1}u^{-2},\,e^{-u/\theta_2}u^{-2}\}$
by a direct substitution. For a multiplicity $k_i\ge2$ the
root of the symbol becomes multiple, and the exact subexponential
structure of the modes (the exponent $\rho$, possible logarithmic
factors and fractional-exponential rank-lowering corrections)
depends on the parameters. It is given by the standard theory of
normal solutions (Thom\'e series) at an irregular singular point
\cite{Wasow1965} (for second-order equations, see
\cite[Chapter~7]{Olver1974}, \cite[Sect.~3.1]{Slavyanov2000}) and is
inessential for our purposes: any such correction is $o(u)$ in the
exponent, so the leading exponential factor $e^{\mu u}$ is
preserved, and the sign of $\mathrm{Re}\,\mu$ --- growth at
$\mu=1/\theta_1$ and decay at $\mu=-1/\theta_2$, the only fact used in
the selection --- remains valid under any multiplicity.
\end{proof}

\begin{remark}\label{rem:exponent-dict}
The exponent $\beta=2a/\sigma^2-1$ here and the exponent
$\gamma=2a_\kappa/\sigma_\kappa^2$ in \cite{PRYu2026} (with $\kappa\in(0,1]$ the
fraction of the reserve invested in the risky asset, so that
$\gamma=2a/\sigma^2$ for $\kappa=1$) are one and the same object in
different parametrizations: $\beta=\gamma-1$, and the power of $u$ in
the tail is the same, $u^{-\beta}=u^{1-\gamma}$, and
$\Phi'(u)\sim D_0u^{-\gamma}$.
\end{remark}

\begin{lemma}[actual fundamental set]\label{lem:actual}
Let $u=\infty$ be an irregular singular point of rank $1$ of a linear
ODE of order $n$ with polynomial coefficients (such are all reduced
equations of Sects.~\ref{sec:reduction}--\ref{sec:asymptotics}), and
let $\mu_1,\dots,\mu_r$ be the roots of its characteristic polynomial
with multiplicities $n_1,\dots,n_r$. Then, in any closed sector with
vertex $u=\infty$ of sufficiently small opening containing the ray
$[u_0,\infty)$, there exists a fundamental set of solutions that splits
into groups of $n_i$ solutions of the form $e^{\mu_iu+o(u)}$. To a
simple root $\mu_i$ there corresponds a solution
$y_i=e^{\mu_iu}u^{\rho_i}(1+o(1))$ whose formal Thom\'e series is its
Poincar\'e asymptotic expansion, which can be differentiated
termwise. A solution carrying a given asymptotic expansion is
fixed up to terms subdominant on the ray (on the positive
half-axis, exponentially small). In particular:
\begin{enumerate}[label=\textup{(\alph*)},leftmargin=2.2em,itemsep=1pt]
  \item To the root $\mu=0$ of the equation for $G=\Phi'$ of
  models~(A) and~(B) (the characteristic polynomials are
  $\mu(1+\theta\mu)^k$ and $\mu(1-\theta\mu)^k$), simple for all $k$,
  there corresponds the solution
  $G_{\mathrm{pow}}(u)=u^{-p}\bigl(\sum_{j<L}a_ju^{-j}+O(u^{-L})\bigr)$
  for each $L\ge1$, where the $a_j$ are the coefficients of the formal
  series of the corresponding model.
  \item In model~(C), equation \eqref{eq:ODE-C-op}, which contains no
  term with $\Phi^{(0)}$, is an equation of order $k_1+k_2+1$ for
  $G=\Phi'$ with the characteristic polynomial $\mu A(\mu)$ and the
  simple root $\mu=0$. To it there corresponds the solution
  $G_{\mathrm{pow}}(u)=u^{-\beta-1}(1+o(1))$, and the algebraic
  branches of Lemma~\ref{lem:fundC} correspond to the exact solution
  $\Phi\equiv1$ and
  \begin{equation}\label{eq:Phibeta}
  \Phi_\beta(u):=\beta\int_u^\infty G_{\mathrm{pow}}(v)\,dv=u^{-\beta}\bigl(1+o(1)\bigr),\qquad \Phi_\beta'=-\beta\,G_{\mathrm{pow}}.
  \end{equation}
\end{enumerate}
Since on $(0,\infty)$ the reduced equations have no singular points,
we may extend the fundamental set to the whole half-line.
\end{lemma}

\begin{proof}
These are the classical asymptotic existence theorems for solutions at
an irregular singular point. For pairwise distinct roots ($k=1$ or
$k_1=k_2=1$), see
\cite[Chapter~V, Sect.~4, Theorem~4.1]{CoddingtonLevinson},
\cite[Theorem~VII-5-3]{HsiehSibuya}, \cite[Sects.~12--15]{Wasow1965}.
Theorem VII-5-3 of \cite{HsiehSibuya} (a Levinson-type theorem for
systems with asymptotically constant coefficients and distinct
eigenvalues) is used for the same purpose in \cite[Sect.~8]{KPuh2020}
and in \cite{PRYu2026}.

We now verify the conditions of these theorems for the reduced
equations. Their coefficients are polynomials in $u$ of degree at most
$2$, the leading coefficient is proportional to $u^2$ and does not
vanish on $(0,\infty)$. After division by the leading coefficient, we write an
equation of order $n$ as the system $Y'=\mathfrak A(u)Y$
for $Y=(y,y',\dots,y^{(n-1)})$ with the matrix
$\mathfrak A(u)=\mathfrak A_0+\mathfrak A_1u^{-1}+\mathfrak A_2u^{-2}$.
Thus the point $u=\infty$ has Poincar\'e rank $1$, and
$\mathfrak A_0$ is the companion matrix of the normalized
characteristic polynomial (for the equation for $\Phi$ it is
$\mathfrak p_2(\mu)$, see below), its eigenvalues being the roots
$\mu_i$ with multiplicities $n_i$. All roots are real: $0$ and
$-1/\theta$ in~(A), $0$ and $1/\theta$ in~(B), $0$, $1/\theta_1$,
$-1/\theta_2$ in~(C). The directions in which two exponentials
$e^{\mu_iu}$, $e^{\mu_ju}$ with $\mu_i\neq\mu_j$ are comparable in
modulus ($\operatorname{Re}((\mu_i-\mu_j)u)=0$) are the rays
$\arg u=\pm\pi/2$, and the positive half-axis lies strictly inside the
sector $|\arg u|<\pi/2$, which contains no such rays. In the case of multiple
roots, the formal system is first split into blocks corresponding to
the distinct roots (the Hukuhara--Turrittin theorem,
\cite[Sect.~19]{Wasow1965}), and the splitting is realized analytically
in a closed sector containing the ray by a transformation
$Y=\Pi(u)Z$ with
$\Pi(u)\sim\Pi_0+\Pi_1u^{-1}+\dots$, $\det\Pi_0\neq0$
\cite[Sects.~12--15, 18--19]{Wasow1965}. The block of the simple root
$\mu=0$ is one-dimensional: $z'=b(u)z$ with
$b(u)\sim\sum_{j\ge1}b_ju^{-j}$, and its solution
$z=u^{b_1}\exp\bigl(\int_\infty^u(b(v)-b_1v^{-1})\,dv\bigr)$ has a
complete asymptotic expansion, termwise differentiable by virtue of
the equation itself. The block of a multiple root $\mu_i$ has the form
$z'=(\mu_iI+B(u))z$, where after a rescaling of the Jordan blocks
$\|B(u)\|\le\eps$ for $u\ge u_0(\eps)$, so that all its solutions satisfy
$|z(u)|=e^{\mu_iu+o(u)}$ (cf.\ the proof of
Proposition~\ref{prop:realization}). This suffices for all assertions
of the lemma and their applications below: the fine structure of the
modes of multiple roots (ramified modes with fractional powers of $u$
in the exponent, whose exponents do not exceed $1-1/n_i$, and
logarithmic factors) is not used below, and the fact that inside the
block of a multiple root the positive half-axis may be a ray of equal moduli for
a pair of conjugate ramified modes does not affect the estimate
$e^{\mu_iu+o(u)}$. For the equations of this paper the branching
exponents are $1-2/n_i$ for $n_i\ge3$, while for $n_i=2$ there is no
branching (modes $e^{\mu_iu}u^{s}$ with exponents $s$ from a quadratic
equation, and with a logarithmic factor if it has a multiple
root). Indeed, the symbol of the equation for
$\Phi$ has the form
$u^2\mathfrak p_2(\mu)+u\mathfrak p_1(\mu)+\mathfrak p_0(\mu)$, where
$A(\mu)$ is the symbol of the annihilator ($(1+\theta\mu)^k$ in~(A),
$(1-\theta\mu)^k$ in~(B), \eqref{eq:symbolC} in~(C)). By the Leibniz
rule $\mathfrak p_2=\frac{1}{2}\sigma^2\mu^2A(\mu)$ and
$\mathfrak p_1=\mu\bigl(aA(\mu)+\sigma^2\mu A'(\mu)\bigr)$, so that at
a root $\mu_i$ of the symbol $A$ of multiplicity $n_i$ the polynomial
$\mathfrak p_1$ has a zero of order $n_i-1$, whereas
$\mathfrak p_0(\mu_i)\neq0$ (for $n_i\ge3$ this is $\alpha$ in the
one-sided models and $\alpha_1(1+\theta_2\mu_i)^{k_2}$ or
$\alpha_2(1-\theta_1\mu_i)^{k_1}$ in model~(C), where
$\alpha_1,\alpha_2>0$ is essential). The balance
$u^2(\mu-\mu_i)^{n_i}\sim \mathfrak p_0(\mu_i)$ gives
$\mu-\mu_i\sim u^{-2/n_i}$, and the term $u\mathfrak p_1$ is
subdominant for $n_i\ge3$. On the positive half-axis, the term with a smaller $\mu_j$ is
maximally subdominant relative to the one with a larger $\mu_i$, that is,
the half-axis is a Stokes line in the sense of Borel summation
(Remark~\ref{rem:realization}). On a ray, a solution
with a given asymptotic expansion is determined only up to terms
subdominant on it: the difference of two solutions with the same
expansion has zero power expansion and is therefore subdominant. The
continuation to $(0,\infty)$ is possible since the leading coefficient
does not vanish there.

For the equations of models~(A) and~(B) the characteristic polynomials
are $\mu(1\pm\theta\mu)^k$, the root $\mu=0$ is simple, and the next
order of the dominant balance gives $\rho=-p$
(Sect.~\ref{subsec:generalk}, and in model~(B) the same computation goes through
with the replacement $\theta\mapsto-\theta$ in the coefficients of the
annihilator). In model~(C), the absence of the term with $\Phi^{(0)}$
was noted after Theorem~\ref{thm:C}. The characteristic polynomial
$\mu^2A(\mu)$ of the equation for $\Phi$ (Lemma~\ref{lem:fundC}),
after passing to $G=\Phi'$, becomes $\mu A(\mu)$ with the simple root
$\mu=0$, and the exponent $-\beta-1$ of the branch $G_{\mathrm{pow}}$
corresponds to the branch $s=-\beta$ for $\Phi$
(Proposition~\ref{prop:exponent}). We integrate the asymptotic
expansion of $G_{\mathrm{pow}}$ over $(u,\infty)$ (legitimate for
$\beta>0$) and obtain \eqref{eq:Phibeta}. The function $\Phi_\beta$
satisfies the equation for $\Phi$, because its derivative satisfies
the equation for $G$.
\end{proof}

\begin{proposition}[realization of the power branch]\label{prop:realization}
Let (A1)--(A2) hold, and, in model~(C), also (A3). Then in each of the
three models the coefficient $c_\beta$ of the power branch
$\Phi_\beta$ in the expansion of the physical solution in the actual
fundamental set of Lemma~\ref{lem:actual} is nonzero, and
\begin{equation}\label{eq:tail-all}
\Psi(u)\sim C_k\,u^{-\beta},\qquad u\to\infty,\qquad C_k=-c_\beta>0,\qquad \beta=\frac{2a}{\sigma^2}-1 .
\end{equation}
In model~(A) we also have, for each $L$,
$\Psi(u)=C_ku^{-\beta}\bigl(\sum_{j<L}\frac{\beta}{\beta+j}\,a_ju^{-j}+O(u^{-L})\bigr)$,
where the $a_j$ are the coefficients of the formal series (for $k=2$,
from \eqref{eq:rec}).
\end{proposition}

\begin{proof}
By Theorems~\ref{thm:A}--\ref{thm:C} the function $G=\Phi'$ satisfies
the reduced equation for $G$ on $(0,\infty)$ and is therefore a linear
combination of the actual fundamental set of Lemma~\ref{lem:actual}.
Integrating, we find that $\Phi$ is a linear combination of a
constant, the branch $\Phi_\beta$ (in models~(A), (B) we define it by
the same formula \eqref{eq:Phibeta} from the corresponding
$G_{\mathrm{pow}}$), and the exponential modes. The growing modes
($e^{u/\theta+o(u)}$ in~(B), $e^{u/\theta_1+o(u)}$ in~(C)) are excluded
by the boundedness of $\Phi$. To see this, we write the equation for
$\Phi$ (of order $K+2$, with the same growing group as that of the
equation for $G$) as the system $Y'=\mathfrak A(u)Y$ for
$Y=(\Phi,\Phi',\dots,\Phi^{(K+1)})$. The analytic block
diagonalization in a sector containing the ray
\cite[Sects.~12--15, 18--19]{Wasow1965}, $Y=\Pi(u)Z$ with
$\Pi(u)=\Pi_0+O(u^{-1})$, $\det\Pi_0\neq0$, separates the growing block
$z'=(\mu I+B(u))z$ with $\mu=1/\theta$ (respectively $1/\theta_1$),
where after a suitable rescaling of the Jordan blocks
$\|B(u)\|\le\eps$ for $u\ge u_0(\eps)$. Then
$\frac{d}{du}|z|^2=2\langle z,(\mu I+B(u))z\rangle\ge2(\mu-\eps)|z|^2$,
whence $|z(u)|\ge e^{(\mu-\eps)(u-u_0)}|z(u_0)|$ for every nonzero
solution of the block, and hence
$|Y(u)|\ge c\,e^{(\mu-\eps)u}$ for every solution $\Phi_{\mathrm g}$
with a nonzero growing component. On the other hand, if such
$\Phi_{\mathrm g}$ were bounded on $[u_0,\infty)$, then so would
its derivatives: the normalized coefficients of the equation are
bounded for $u\ge u_0$, and by the interpolation inequality on the
intervals $[x,x+1]$ (here $\|\cdot\|$ is the uniform norm on such an
interval, and the constants do not depend on $x$),
$\|\Phi_{\mathrm g}^{(j)}\|\le\eta\|\Phi_{\mathrm g}^{(K+2)}\|+C_\eta\|\Phi_{\mathrm g}\|$,
$0<j<K+2$, so that, together with the equation (all the interval norms
are finite),
$\|\Phi_{\mathrm g}^{(K+2)}\|\le C\|\Phi_{\mathrm g}\|$ with a constant
independent of $x$ --- in contradiction with the growth of $|Y|$. Thus
the growing component of a bounded solution vanishes, whatever the
logarithmic factors and the complex exponents inside the block
may be, and the bounded solutions are exactly the linear combinations
of a constant, the branch $\Phi_\beta$, and the decaying modes. The
decaying modes have the form $e^{-u/\theta+o(u)}$ (in~(C):
$e^{-u/\theta_2+o(u)}$) and are therefore $O(e^{-\delta u})$ for some
$\delta>0$. In all three models this yields
\[
\Phi=c_0+c_\beta\Phi_\beta+O(e^{-\delta u}),
\]
and in model~(B) there are no decaying modes and the remainder is
absent. From $\Phi(\infty)=1$ it follows that $c_0=1$, i.e.,
$\Psi=-c_\beta\Phi_\beta+O(e^{-\delta u})$. We still need to prove that
$c_\beta\neq0$.

\emph{Model~(B).} The solutions bounded as $u\to\infty$ form
$\mathrm{span}\{1,\Phi_\beta\}$, and $\Phi=1+c_\beta\Phi_\beta$. If
$c_\beta=0$, then $\Phi\equiv1$, which contradicts $\Phi(0+)=0$ (and
$\Psi>0$). Thus $c_\beta\neq0$.

\emph{Models~(A), (C).} If $c_\beta=0$, then
$\Psi=O(e^{-\delta u})$. This contradicts the power-law lower bound
\cite[Proposition~4.1]{KPuh2020}: if $c<0$ or the downward jumps are
unbounded, then there exists a $\beta_*>0$ such that
$\liminf_{u\to\infty}u^{\beta_*}\Psi(u)>0$. In model~(C) the
conditions of that proposition are satisfied for any $c$, including
$c=0$: the risk process is a compound Poisson process with drift and
jumps of both signs, as in \cite{KPuh2020}, and the downward jumps
$\xi^{(2)}\sim\Gamma(k_2,\theta_2)$ are unbounded. Model~(A) does not
satisfy the standard assumption of \cite{KPuh2020} about jumps of
both signs, but the proof of Proposition~4.1 of \cite{KPuh2020} does
not use upward jumps: it relies only on the fact that a jump
together with the drift contribution before it is not bounded from
below (in model~(A) this is provided by the unbounded claims
$\xi\sim\Gamma(k,\theta)$), and on the Brownian factor between the
jumps. Therefore it carries over verbatim to model~(A). For jumps
with a finite moment of order $\beta$ we also obtain, from
\cite[Theorem~1.1]{KP2020},
$0<\liminf_{u\to\infty}u^{\beta}\Psi(u)\le\limsup_{u\to\infty}u^{\beta}\Psi(u)<\infty$.
Power-law lower bounds were also obtained in
\cite{KalashnikovNorberg2002} (general models with investments) and,
for model~(A), in \cite{FKP2002} (exponential jumps) and
\cite{PergamenshchikovZeitouni2006}.

\emph{The sign and the expansion.} From $\Psi\ge0$ and $\Phi_\beta>0$
for large $u$ it follows that $-c_\beta\ge0$, hence
$C_k:=-c_\beta>0$ in all three models, which gives
\eqref{eq:tail-all}. In model~(A), the termwise integration of the
asymptotic expansion of $G_{\mathrm{pow}}$ over $(u,\infty)$,
legitimate for power-law asymptotics with exponents $p+j>1$, gives
$\Phi_\beta=\beta\sum_{j<L}\frac{a_j}{\beta+j}u^{-\beta-j}+O(u^{-\beta-L})$,
whence the complete expansion follows, and the amplitude $\varpi$ of
Theorem~\ref{thm:formal-k2} equals $\beta C_k$.
\end{proof}

For one-sided models the asymptotics $\Psi(u)\sim Cu^{-\beta}$,
$C>0$, is known for general claim distributions: in model~(A) under
$\bbe[\xi^{\beta}]<\infty$ \cite[Theorem~1]{P05143}, in model~(B)
under $\bbe[\xi]<\infty$ \cite[Theorem~1]{P2601}, and, under
$\bbe[\xi^{\beta}]<\infty$, \cite[Theorem~1.1]{KP2020}.
Proposition~\ref{prop:realization} gives for them (for
Erlang-distributed claims) a different proof of the asymptotic part.
The smoothness, the value of $\Phi(0+)$, and the lower bound are taken
from \cite{KPuh2020,P2601,P05143}. For the mixed model the asymptotics
with $C>0$ are also known: they follow from
\cite[Theorem~2.1]{BKM2011} (the Cram\'er-type condition on the
logarithm of the price of the risky asset with root $\beta$, the
moment conditions, and the absolute continuity of its law are
satisfied, see also \cite{Nyrhinen2001} and, for independent L\'evy
processes, \cite[Theorem~3.2]{Paulsen2002}), whereas for a general
L\'evy return process with jumps of both signs, only two-sided
estimates of order $u^{-\beta}$ were established in
\cite[Theorem~1.1]{KP2020}. Analytically it was obtained for
$k_1=k_2=1$ \cite{KPuh2020} and, for rational Laplace transforms, for
$\beta\in(0,1)$, $c\neq0$ \cite{Antipov2026}.
Proposition~\ref{prop:realization} gives an analytical proof for all
$k_1,k_2$, $\beta>0$, $c\in\bbr$, in which, apart from smoothness and
the boundary condition $\Phi(\infty)=1$ (Sect.~\ref{sec:model}), only
the power-law lower bound \cite[Proposition~4.1]{KPuh2020} is taken
from probabilistic results (the nonvanishing of $c_\beta$). The new
element is the complete expansion up to a constant factor $C$ in all
three models (in model~(A) with $k=1$ it was obtained in
\cite{PRYu2026}).

\begin{theorem}[the physical solution and its tail, model~(C)]\label{thm:C-asymp}
Under assumptions (A1)--(A3) the physical survival function $\Phi$ of
model~(C) is singled out in terms of the actual fundamental set
(Lemmas~\ref{lem:fundC}, \ref{lem:actual}) as follows: the
coefficients of all $k_1$ growing branches $e^{u/\theta_1+o(u)}$
vanish (the boundedness of $\Phi$), the coefficient of the constant
branch is $1$ (the condition $\Phi(\infty)=1$), and the coefficient
$c_\beta$ of the power branch is nonzero
(Proposition~\ref{prop:realization}), so that
\begin{equation}\label{eq:C-tail}
\Psi(u)\sim C\,u^{-\beta},\qquad u\to\infty,\quad C=-c_\beta>0,\quad \beta=\frac{2a}{\sigma^2}-1,
\end{equation}
with the same exponent as in the one-sided models, while the
contribution of the $k_2$ decaying branches is
$O(e^{-(1/\theta_2-\eps)u})$ for any $\eps>0$.
\end{theorem}

\begin{proof}
The growing branches $e^{u/\theta_1+o(u)}$ are excluded by the
boundedness of $\Phi\in[0,1]$ --- this is the feature distinguishing
models~(B) and~(C) from model~(A): the upward-jump operator
$\cI_{k_1}$ generates $k_1$ growing modes, and all of them are cut
off. After their exclusion,
$\Phi=c_0+c_\beta\Phi_\beta+O(e^{-(1/\theta_2-\eps)u})$ (for a
multiplicity $k_2\ge2$ the decaying modes contain ramified
corrections $e^{o(u)}$, hence the $\eps$ in the remainder). The
condition $\Phi(\infty)=1$ gives $c_0=1$, and then
$\Psi=1-\Phi=-c_\beta u^{-\beta}(1+o(1))+O(e^{-(1/\theta_2-\eps)u})$.
From $\Psi\ge0$ and $\Phi_\beta>0$ for large $u$ it follows that
$c_\beta\le0$, while $c_\beta\neq0$ is given by
Proposition~\ref{prop:realization}, which yields \eqref{eq:C-tail}. The
reduced ODE expresses only a necessary condition
(Remark~\ref{rem:selection}): its $(k_1+k_2+2)$-dimensional solution
space is wider than that of the original IDE, and the data at infinity
(the boundedness, $\Phi(\infty)=1$) fix only $k_1+1$ degrees of
freedom. The IDE itself determines which of the $k_2$ subdominant modes enter
the physical $\Phi$ and what the amplitude $c_\beta$ is: this is fixed
by at most $k_2$ effective conditions \eqref{eq:converse} of
Proposition~\ref{prop:converse} (Remark~\ref{rem:converse}) and by the
behavior at $u=0$: for $c<0$ by the condition $\Phi(0+)=0$, for
$c>0$ by boundedness, which excludes the mode \eqref{eq:C-irr}
(Proposition~\ref{prop:C-zero}), and for $c=0$ by boundedness at zero
(Remark~\ref{rem:c0}).
\end{proof}

As in model~(A), we obtain a complete asymptotic expansion up to a
constant factor $C$. The root $\mu=0$ is simple, so the formal Thom\'e
series $G_{\mathrm{pow}}=u^{-\beta-1}\sum_{j\ge0}a_ju^{-j}$ with
$a_0=1$ and coefficients uniquely computed by the recurrence of the
equation for $G=\Phi'$ (in normalized form the factor multiplying $a_j$ in it
equals $-j/\bigl((- \theta_1)^{k_1}\theta_2^{k_2}\bigr)\neq0$ for
$j\ge1$, as in Proposition~\ref{prop:general-k}) is its asymptotic
expansion (Lemma~\ref{lem:actual}), and the termwise integration, as
in the proof of Proposition~\ref{prop:realization}, gives
$\Psi(u)=Cu^{-\beta}\bigl(1+\sum_{1\le j<L}\frac{\beta a_j}{\beta+j}u^{-j}\bigr)+O(u^{-\beta-L})$
for any $L\ge1$. The same is true in model~(B), where
$\Psi=-c_\beta\Phi_\beta$ exactly.

\subsection{The singular point $u=0$ in the mixed model}\label{subsec:zero}

\begin{proposition}[behavior at the singular point $u=0$, model~(C)]\label{prop:C-zero}
Let $c\neq0$. The point $u=0$ is an irregular singular point of rank~1
of the reduced ODE \eqref{eq:ODE-C-op}: the leading coefficient is
proportional to $u^2$, while the coefficient of $\Phi^{(M-1)}$,
$M:=k_1+k_2+2$, is nonzero at zero. Formal power solutions
$u^\nu\sum_{n\ge0}b_nu^n$ exist for $\nu=0,1,\dots,k_1+k_2$ --- the
lowest-order condition is the indicial equation
\begin{equation}\label{eq:C-ind0}
\nu(\nu-1)\cdots\bigl(\nu-(k_1+k_2)\bigr)=0,
\end{equation}
and they form a $(k_1+k_2+1)$-parameter family. The convergent series form
a subspace of it whose dimension is at least $K=k_1+k_2$ and which
contains the constant. The dimension equals $K$ exactly when at least
one formal series diverges. In that case (the factorially growing solution
of the recurrence gives $|b_{n+1}/b_n|\sim n\sigma^2/(2|c|)$, a Gevrey
class of order $1$) the divergent tails of all formal solutions are
proportional. In either case, since $K=k_1+k_2\ge2$, there exist
nonconstant solutions analytic at zero (and entire, because the
equation has no other finite singular points). The solution
corresponding to the formal solution of exponential type has an
essential singularity at zero,
\begin{equation}\label{eq:C-irr}
\Phi_{\mathrm{irr}}(u)=u^{1-\beta}\exp\Bigl(\frac{2c}{\sigma^2u}\Bigr)\bigl(1+o(1)\bigr),\qquad u\to0+.
\end{equation}
For $c>0$ it grows without bound and is excluded by the boundedness of
$\Phi$, while for $c<0$ it is exponentially small together with all
its derivatives and is not ruled out by any regularity condition. The
limit $u\to0+$ in the IDE \eqref{eq:IDE-C} gives the boundary relation
\begin{equation}\label{eq:C-bc0}
c\,\Phi'(0+)=(\alpha_1+\alpha_2)\Phi(0+)-\alpha_1\,\cI_{k_1}\Phi(0).
\end{equation}
For $c<0$ we have $\Phi(0+)=0$ (property (iii) of
Sect.~\ref{sec:model}), and then \eqref{eq:C-bc0} determines
$\Phi'(0+)=\alpha_1\cI_{k_1}\Phi(0)/|c|>0$. For $c>0$, from
\eqref{eq:C-bc0}, $\Phi'\ge0$ and $\alpha_1>0$ it follows that
$\Phi(0+)\ge\alpha_1\cI_{k_1}\Phi(0)/(\alpha_1+\alpha_2)>0$.
\end{proposition}
\begin{proof}
The leading coefficient of the ODE \eqref{eq:ODE-C-op} is
$(-1)^{k_1}\frac{1}{2}\sigma^2\theta_1^{k_1}\theta_2^{k_2}u^2$, and
the coefficient of $\Phi^{(M-1)}$ at zero is
$(-1)^{k_1}c\,\theta_1^{k_1}\theta_2^{k_2}$ (the contribution of the
drift $(au+c)\Phi'$ through the leading term of the annihilator). They
share the sign factor. The normalized coefficient of $\Phi^{(M-1)}$
has a second-order pole at zero. Hence, by the Fuchs criterion, the
point is irregular, and after the substitution $t=1/u$ the exponential
factor is linear in $t$, i.e., the rank is $1$. Substituting
$\Phi\sim u^\nu$, we find the lowest order $u^{\nu-(M-1)}$ given by
a single term $c\Phi^{(M-1)}$, which yields \eqref{eq:C-ind0}. For
$\nu\in\{0,\dots,k_1+k_2\}$ the coefficients $b_0,\dots,b_{k_1+k_2}$
are free, while for $n>k_1+k_2$ the principal part of the recurrence
has the form
$c\,n^{\underline{M-1}}\,b_n=-\frac{1}{2}\sigma^2(n-1)^{\underline{M}}\,b_{n-1}+\dots$
($x^{\underline m}$ is the falling factorial), whence
$b_n/b_{n-1}\sim-n\sigma^2/(2c)$ for solutions of the recurrence that
contain the factorially growing component. The constant $\Phi\equiv1$ is an
exact solution, because the term in $\Phi^{(0)}$ is absent. The
dimension of the subspace of convergent series is given by the
Malgrange index theorem \cite{Malgrange1974}: if
$\mathcal T=\sum_{i=0}^{M}t_i(u)D^i$ is an operator with coefficients
holomorphic at zero, $t_M\not\equiv0$, and $\operatorname{ord}t_i$ is
the order of the zero of $t_i$ at $u=0$ ($\operatorname{ord}0=+\infty$),
then $\mathcal T$ is Fredholm both on $\bbc[[u]]$ and on the convergent
power series $\bbc\{u\}$, with the indices
\[
\chi(\mathcal T,\bbc[[u]])=\max_{0\le i\le M}\bigl(i-\operatorname{ord}t_i\bigr),\qquad \chi(\mathcal T,\bbc\{u\})=M-\operatorname{ord}t_M,
\]
and their difference equals the irregularity of the point $u=0$. For
the operator of equation \eqref{eq:ODE-C-op} we have
$\operatorname{ord}t_M=2$, $\operatorname{ord}t_{M-1}=0$ (the
coefficient $(-1)^{k_1}c\,\theta_1^{k_1}\theta_2^{k_2}\neq0$ for
$c\neq0$) and $\operatorname{ord}t_i\ge0$ for $i<M-1$, so that
$\chi(\bbc[[u]])=M-1=k_1+k_2+1$,
$\chi(\bbc\{u\})=M-2=k_1+k_2$, and the irregularity equals $1$ (the
unique exponential mode \eqref{eq:C-irr} of rank $1$, and the Newton
polygon has one sloping edge of length $1$ and slope $1$). This agrees
with a direct count: on $\bbc[[u]]$ we find the operator surjective
(the coefficient $c\,n^{\underline{M-1}}$ of the leading unknown $b_n$
is nonzero for $n\ge M-1$), and its kernel has dimension
$k_1+k_2+1$.
The dimension of the kernel on $\bbc\{u\}$ is therefore at least the
index $k_1+k_2$. It equals $k_1+k_2$ exactly when the kernel on
$\bbc\{u\}$ is strictly smaller than the kernel on $\bbc[[u]]$,
that is, when at least one formal series diverges. In this case the divergent
tails of all formal solutions are proportional to a single factorially
growing tail, and the divergence is governed by one linear functional
of the free parameters $b_0,\dots,b_{k_1+k_2}$ (the Stokes
multiplier). An index count alone cannot rule out the vanishing of the
Stokes multiplier: when the irregularity equals $1$ all formal series may converge
(as for the operator $(u^2D+1)D$). For \eqref{eq:ODE-C-op} this would
mean the convergence of all $k_1+k_2+1$ formal series. At none of the
tested parameter values is this observed (for $k=1$ in model~(A) the
divergence for generic parameters was noted in \cite{PRYu2026}).
However, we do not have a proof for all parameter values (for fixed
$k_1,k_2$ and the sign of $c$, the set of parameters at which all the
series converge is the zero set of the Stokes multiplier, apparently
analytic in the parameters, and the fact that it is a proper subset of this component would
follow from a single rigorously established example of divergence in
it).

For the solution with an essential singularity we substitute
$\Phi=e^{S}$ with $S=\lambda/u+m\ln u$. The two leading terms
$\frac{1}{2}\sigma^2u^2\Phi^{(M)}+c\Phi^{(M-1)}$ give, in the most
singular order $u^{2-2M}$, a factor proportional to
$\lambda^{M-1}(2c-\lambda\sigma^2)$, whence $\lambda=2c/\sigma^2$.
Equivalently, for $y=\Phi^{(M-1)}$ the balance
$\frac{1}{2}\sigma^2u^2y'+cy=0$ integrates exactly:
$y=\exp(2c/(\sigma^2u))$. The next order $u^{3-2M}$ receives contributions from
the term $m\,u^{-1}$ of $S'$, the terms
$\binom{M}{2}S'^{M-2}S''$ and $\binom{M-1}{2}S'^{M-3}S''$ of the
expansions of $\Phi^{(M)}/\Phi$ and $\Phi^{(M-1)}/\Phi$, the term
$(M-2)\sigma^2u\,\Phi^{(M-1)}$ from the Leibniz rule, and the drift
$au\,\Phi^{(M-1)}$. Their sum vanishes at
$m=2-2a/\sigma^2=1-\beta$, independently of $k_1$, $k_2$, which gives
\eqref{eq:C-irr}. The relation \eqref{eq:C-bc0} is the limit of
\eqref{eq:IDE-C} as $u\to0+$: by property (iv) of Sect.~\ref{sec:model},
the derivatives $\Phi'$, $\Phi''$ are bounded, so that
$\frac{1}{2}\sigma^2u^2\Phi''\to0$ and $au\Phi'\to0$, while
$\cJ_{k_2}\Phi(u)\to0$. The limit $\Phi'(0+)$ exists because $\Phi''$
is bounded (property (iv), $c\neq0$), so $\Phi'$ is Lipschitz on
$(0,1)$. For $c>0$ the growing mode \eqref{eq:C-irr} is excluded by
boundedness. The assertions about $\Phi(0+)$ follow from
\eqref{eq:C-bc0} immediately: for $c>0$ and $\Phi(0+)=0$ we would have
$c\Phi'(0+)=-\alpha_1\cI_{k_1}\Phi(0)<0$ contrary to $\Phi'\ge0$.
\end{proof}

\begin{remark}\label{rem:c0}
For $c=0$ the coefficient of $\Phi^{(M-1)}$ vanishes at $u=0$, the
Fuchs condition holds, and the point $u=0$ becomes a regular singular
one: all solutions are convergent Frobenius series (possibly with
logarithms), and the mode \eqref{eq:C-irr} is absent. The indicial
equation has the roots $0,1,\dots,k_1+k_2-1$ and two roots
$\nu_-<0<\nu_+$ of the quadratic trinomial
$\sigma^2\nu^2+(2a-\sigma^2)\nu-2(\alpha_1+\alpha_2)$. The root
$\nu_-$ corresponds to the branch $u^{\nu_-}$, unbounded at the origin
and therefore excluded. The physical solution is bounded and is
therefore a combination of the Frobenius branches with exponents
$0,1,\dots,k_1+k_2-1$ and $\nu_+$:
$\Phi(u)=\Phi(0+)+O(u^{\min(1,\nu_+)}|\ln u|)$. After termwise
differentiation of the Frobenius expansions, $u\Phi'(u)\to0$ and
$u^2\Phi''(u)\to0$ as $u\to0+$. The limit in \eqref{eq:IDE-C} gives
$(\alpha_1+\alpha_2)\Phi(0+)=\alpha_1\cI_{k_1}\Phi(0)$, i.e.,
$\Phi(0+)=\alpha_1\cI_{k_1}\Phi(0)/(\alpha_1+\alpha_2)>0$: survival
from zero initial capital is possible only if the first jump is
directed upward. The first part of property (iv), the boundedness of
the derivatives, does not apply here: for $\nu_+<1$, i.e., for
$a>\alpha_1+\alpha_2$, the derivatives $\Phi'$, $\Phi''$ are unbounded
near zero. For $c=0$ the IDE \eqref{eq:IDE-C} is then an
Euler-type equation
$\frac{1}{2}\sigma^2u^2\Phi''+au\Phi'-\bar\alpha\Phi=-h$ with
$\bar\alpha:=\alpha_1+\alpha_2$ and
$h:=\alpha_1\cI_{k_1}\Phi+\alpha_2\cJ_{k_2}\Phi$. Here $h$ is bounded
and nondecreasing, as a sum of convolutions of a nondecreasing
function $\Phi$. For $u>0$ the inequality $h(u)-h(0)\ge\alpha_2\cJ_{k_2}\Phi(u)>0$ holds
and, by Lemma~\ref{lem:base},
$h'=\alpha_1(I_{k_1}-I_{k_1-1})/\theta_1+\alpha_2(J_{k_2-1}-J_{k_2})/\theta_2$,
so that $|h'|\le\alpha_1/\theta_1+\alpha_2/\theta_2$ since
$0\le\Phi\le1$, and $h$ is Lipschitz (the boundedness of $\Phi'$ is
not needed for this). The unique bounded solution of this Euler-type
equation is given by
\[
\Phi(u)=\frac{2}{\sigma^2(\nu_+-\nu_-)}\Bigl[u^{\nu_-}\!\int_0^u x^{-\nu_--1}h(x)\,dx+u^{\nu_+}\!\int_u^\infty x^{-\nu_+-1}h(x)\,dx\Bigr],
\]
whence, taking into account $-\nu_-\nu_+=2\bar\alpha/\sigma^2$,
$\Phi(u)=h(0)/\bar\alpha+\frac{2H}{\sigma^2(\nu_+-\nu_-)}\,u^{\nu_+}+O(u)$,
where $H=\int_0^\infty x^{-\nu_+-1}\bigl(h(x)-h(0)\bigr)\,dx>0$ for
$\nu_+<1$. This expansion also gives the value
$\Phi(0+)=h(0)/\bar\alpha$. For $\nu_+=1$ the growth of $\Phi'$ is logarithmic, and for
$1<\nu_+<2$, $\Phi'$ stays bounded, while $\Phi''$, in general, does
not.
\end{remark}

\begin{remark}\label{rem:C-open}
In Theorem~\ref{thm:C-asymp}, after excluding the growing modes and
normalizing $c_0=1$, there remain $k_2+1$ unknowns: the amplitude
$c_\beta$ and the coefficients of the $k_2$ decaying modes, and they
are not determined by the reduced ODE. They are fixed by at most
$k_2$ effective conditions \eqref{eq:converse}
(Proposition~\ref{prop:converse}, Remark~\ref{rem:converse}), i.e.,
by a return to the IDE, and by one condition at zero:
$\Phi(0+)=0$ for $c<0$, boundedness for $c>0$ (the exclusion of the
mode \eqref{eq:C-irr}) and for $c=0$ (the exclusion of the branch
$u^{\nu_-}$, Remark~\ref{rem:c0}). For $k_1=k_2=1$ (exponential
jumps of both signs) a complete description of the solution space
and the asymptotics with $C>0$ are given in \cite[Sects.~7--8]{KPuh2020}.
The exact solution through Heun functions \cite{PRYu2026} refers to
the one-sided model~(A) with exponentially distributed claims. For
general $k_1,k_2$, an explicit solution of the connection problem
between $u=0$ and $u=\infty$ --- the computation of the coefficients
$c_\beta$ and $C_k$ themselves --- remains open. This affects neither
the exponent $\beta$ nor the positivity of $C_k$
(Proposition~\ref{prop:realization}).
\end{remark}

\section{Test cases}\label{sec:examples}

For $k=1$: $\cN(f)=f+\theta f'$ and
$J_1(u)=\theta^{-1}\int_0^u\Phi(u-y)e^{-y/\theta}\,dy$. The formula
of Theorem~\ref{thm:A} for $k=1$ gives an ODE of order $3$ for $\Phi$
or, equivalently, of order $2$ for $G=\Phi'$:
\begin{equation}\label{eq:ODE-k1}
\frac{1}{2}\sigma^2\theta u^2 G''+\Bigl(\frac{1}{2}\sigma^2 u^2+\theta(\sigma^2 u+au+c)\Bigr)G'+\bigl(\theta a+au+c-\theta\alpha\bigr)G=0,
\end{equation}
which coincides with equation (13) of \cite{FKP2002} and admits an
exact solution through Heun functions, as shown in \cite{PRYu2026}.
We use this as a consistency check of the formulas of
Theorem~\ref{thm:A} in the simplest case.

For $k=2$ the operator is $\cN(f)=f+2\theta f'+\theta^2 f''$. By
Theorem~\ref{thm:A}:
\begin{align*}
g_4(u)&=\frac{1}{2}\sigma^2\theta^2 u^2,\\
g_3(u)&=\sigma^2\theta u^2+\theta^2(2\sigma^2 u+au+c),\\
g_2(u)&=\frac{1}{2}\sigma^2 u^2+2\theta(\sigma^2 u+au+c)+\theta^2(\sigma^2+2a-\alpha),\\
g_1(u)&=2\theta a+au+c-2\theta\alpha,
\end{align*}
which agrees with the coefficients of \eqref{eq:ODE-k2}. By
Theorem~\ref{thm:formal-k2} and Proposition~\ref{prop:realization},
$\Psi(u)\sim C_2u^{-\beta}$ with $C_2>0$, and we compute the coefficients of the
complete expansion from the recurrence \eqref{eq:rec}:
$a_1=2(c-2\alpha\theta)/\sigma^2$,
$a_2=-2(3a\alpha\theta^2-4\alpha^2\theta^2+4\alpha c\theta-c^2)/\sigma^4$
(Remark~\ref{rem:dominant}).

For $k_1=k_2=1$ the operator form \eqref{eq:ODE-C-op} gives an ODE of
order $4$ for $\Phi$ (of order $3$ for $G$), which for
$\theta_1=\theta_2$ coincides with the expanded form
\eqref{eq:ODE-C}. This case agrees with equation (7.2) for $G=\Phi'$
of the mixed model in \cite{KPuh2020} for exponential jumps with
arbitrary $\theta_1$, $\theta_2$ (up to an overall sign and a
renumbering of the indices: there, index~1 corresponds to downward
jumps), which provides a consistency check for Theorem~\ref{thm:C}.

\section{Discussion}\label{sec:discussion}

In this paper we carried out the reduction of the IDE to a linear ODE
as a necessary consequence for all three variants of the risk model
with investments and Erlang-distributed jumps. The main results:

\begin{enumerate}[label=(\roman*)]
  \item We constructed annihilating operators $\cN$ and $\wt\cN$ of
  order $k$ that explicitly eliminate the nonlocal integral in the
  IDE.
  \item The coefficients of the resulting ODEs are given in closed
  form in terms of binomial coefficients for the one-sided models and
  for the mixed model with $k_1=k_2$ (in the general case, through
  the composition of operators). For $k_1=k_2$ the equation of the
  mixed model coincides with equation (6) of Antipov
  \cite{Antipov2026} for densities with rational Laplace transform.
  The remaining equations follow from the same construction.
  \item We established the exponent of the power branch
  $\beta=2a/\sigma^2-1$, independent of $k$, and proved the
  realization of the tail $\Psi(u)\sim C_ku^{-\beta}$, $C_k>0$, for
  all $k$ and all $\beta>0$ (Proposition~\ref{prop:realization}). For
  the one-sided models this gives a different proof of the asymptotic part
  of the results \cite{KP2020,P2601,P05143} (in the annuity model,
  without any external tail bound, from the ODE and the condition
  $\Phi(0+)=0$, and in non-life insurance with the power-law lower
  bound \cite[Proposition~4.1]{KPuh2020}). For the mixed model it gives
  an analytical proof, with the same lower bound, of a result known
  from the probabilistic approach
  \cite{BKM2011,Nyrhinen2001,Paulsen2002}, extending \cite{KPuh2020}
  to Erlang-distributed jumps and \cite{Antipov2026} to $\beta\ge1$
  and $c=0$. In all three models we obtained a complete asymptotic
  expansion up to a multiplicative constant $C_k$ (for
  non-life insurance it was known for $k=1$ from \cite{PRYu2026}, and for $k=2$
  we give its coefficients by an explicit recurrence).
  \item For the mixed model~(C) we classified the fundamental set at
  $u=\infty$. The boundedness of the physical solution excludes the
  $k_1$ growing modes, and the remaining connection coefficients are
  fixed by at most $k_2$ effective conditions of
  Proposition~\ref{prop:converse} (Remark~\ref{rem:converse}) and one
  condition at zero: $\Phi(0+)=0$ for $c<0$, boundedness for $c>0$
  (the exclusion of the mode \eqref{eq:C-irr}) and for $c=0$ (the
  exclusion of the unbounded Frobenius branch).
  \item We obtained a converse criterion for bounded solutions:
  IDE $\Leftrightarrow$ ODE together with the vanishing of the
  residual of the IDE and its derivatives up to order $K-1$ at zero
  (here $K$ is the order of the annihilator, and the number of effective
  conditions is at most $k$, zero, and at most $k_2$ in models~(A),
  (B), (C) respectively, Remark~\ref{rem:converse}), generalizing the
  relation $c\Phi'(0+)=\alpha\Phi(0+)$ of \cite{PRYu2026}.
\end{enumerate}

Natural directions for further work are: (a) explicit solutions
(analogues of the Heun functions) and the computation of the
constants $C_k$ for $k\ge 2$, (b) the asymptotics as $u\to0+$ and
the summation of divergent series at zero (a Gevrey class of order
$1$, cf.\ \cite{PRYu2026}), (c) weakening the smoothness conditions
(A2) in the spirit of \cite{P05143} for distributions with a less
regular density, and (d) for the mixed model~(C), an explicit solution of
the connection problem between $u=0$ and $u=\infty$
(Remark~\ref{rem:C-open}) and the structure of ramified modes at
multiple roots of the symbol.

\vspace{6pt}

\section*{Declarations}

\paragraph{Conflict of interest.} The author declares no conflict of interest.

\paragraph{Availability of data and materials.} No datasets were generated or analyzed during the current study.


\end{document}